%% file: main.tex
\documentclass[10pt]{scrartcl} 

\usepackage[utf8]{inputenc}
\usepackage[T1]{fontenc}
\usepackage{lmodern}

\usepackage{amsfonts,amssymb,amsmath,amsthm}
\usepackage{graphicx}  
   \usepackage[dvipsnames]{xcolor}
\usepackage{natbib}
     \usepackage{authblk}
\usepackage[]{enumitem}
\usepackage{tabularx}
\usepackage{subfig}
\usepackage{mathtools}
\usepackage{url}
\usepackage{algorithm} 
\usepackage{algpseudocode} 
\usepackage{comment}
\usepackage{pgfplotstable}
\usepackage{csvsimple}
\usepackage{array,multirow,makecell}

\newcommand{\eps}{{\ensuremath{\varepsilon}}}

\DeclareMathOperator*{\cl}{cl}

\newcommand{\WS}{\mathit{WS}} 
\newcommand{\AWT}{AWT}            
\newcommand{\LBS}{\mathcal{L}}
\newcommand{\UBS}{\mathcal{U}}

\newcommand{\R}{\mathbb{R}}

\newcolumntype{C}[1]{>{\centering\arraybackslash }b{#1}}

\usepackage{tikz, pgfplots}
\usetikzlibrary{arrows, shapes.misc, chains, patterns}
\usetikzlibrary{decorations.pathreplacing}
\usepackage{xifthen}

\tikzset{
  candidat/.style={rectangle, inner sep=0pt, minimum size=0.1cm, draw=gray, fill=gray},
  nds/.style={circle, inner sep=0pt, minimum size=0.12cm, draw=black, fill=black},
  ndns/.style={thick, draw=red, cross out, inner sep=0pt, minimum width=4pt, minimum height=4pt},
 ndns/.style={rectangle, inner sep=0pt, minimum size=0.12cm, draw=black, fill=black},
 test/.style={circle, inner sep=0pt, minimum size=0.12cm, draw=black, fill=black},
}
\pgfplotsset{compat=1.8}
\definecolor{ChromeYellow}{rgb}{1.0, 0.65, 0.0}

\begin{document}

\title{Augmenting Bi-objective Branch and Bound by Scalarization-Based Information}

\author{Julius Bau{\ss}\thanks{Corresponding author} }
\author{Michael Stiglmayr}
\affil{University of Wuppertal, School of Mathematics and Natural Sciences, IMACM, Gaußstr.~20, 42119 Wuppertal, Germany\\
\{bauss,stiglmayr\}@math.uni-wuppertal.de
}

\date{}
\maketitle
\begin{abstract}
\input{abstract}

\end{abstract}

\section{Introduction}\label{sec:intro}
\input{sec_intro}

\section{Preliminaries}\label{sec:def}
\input{sec_def}

\section{A Generic Multi-objective Branch and Bound Framework}\label{sec:BB}
\input{sec_BB}

\section{Using Objective Space Information in Multi-objective Branch and Bound}\label{sec:hyb}
\input{sec_hyb}

\section{Numerical Results}\label{sec:res}
\input{sec_res}

\section{Conclusion and Outlook}\label{sec:dis}
\input{sec_concl}

\section*{Acknowledgment}
The authors thankfully acknowledge financial support by Deutsche Forschungsgemeinschaft, project number~KL~1076/11-1.


\bibliographystyle{abbrvnat}
\bibliography{bib}


\end{document}

%% file: abstract.tex
While branch and bound based algorithms are a standard approach to solve single-objective (mixed-)integer optimization problems, multi-objective branch and bound methods are only rarely applied compared to the predominant objective space methods. 
In this paper we propose modifications to increase the performance of multi-objective branch and bound algorithms by utilizing scalarization-based information. We use the hypervolume indicator as a measure for the gap between lower and upper bound set to implement a multi-objective best-first strategy. By adaptively solving scalarizations in the root node to integer optimality we improve both, upper and lower bound set. The obtained lower bound can then be integrated into the lower bounds of all active nodes, while the determined solution is added to the upper bound set. 
Numerical experiments show that the number of investigated nodes can be significantly reduced by up to 83\% and the total computation time can be reduced by up to 80\%.

%% file: sec_intro.tex
Many optimization problems occurring in real-word applications include a conflict of interests and goals, or secondary objectives, in a word, they are multi-objective. Thus, there is (in general) not one solution that optimizes all objectives at once. Following the \textit{a posteriori} paradigm of decision making, we aim at determining the set of so-called \emph{efficient solutions} or the images, the so-called \emph{non-dominated points}, which cannot be improved in one objective without deterioration in at least one other objective. Thus, efficient solutions are reasonable choices for decision makers.

As we are considering specifically bi-objective integer linear programs and their solution with multi-objective branch and bound methods, the following literature survey will also focus on this and closely related topics. A comprehensive introduction to multi-objective optimization in general is given, e.g., in \citet{steuer86multiple,ehrgott05multicriteria}.

Solution approaches for multi-objective optimization problems are often categorized in: \emph{objective space} and \emph{decision space methods}. 
Objective space methods scalarize the underlying problem, i.\,e., it is replaced by a series of single-objective problems to determine successively the set of efficient solutions. In the case of multi-objective integer programming, these scalarized problems can be solved with commercial integer programming solvers like CPLEX or Gurobi. The utilization of these optimized, single-criteria solvers are a major advantage and one of the reasons why those methods are predominant in multi-objective optimization.

There are numerous objective space methods and a popular one is the $\eps$-constraint method that was introduced for two objectives by \citet{Haimes1971on}. In every iteration the first objective is optimized with an updated constraint to ensure an improvement regarding the second objective. In \citet{Laumanns2006an} an extension to three and more objectives is presented. Many approaches based on the $\eps$-constraint method have been published in the last decades, for example \citet{Boland2017the} and \citet{Kirlik2014a} combine the method with reduction of dimension in the tri- respectively multi-dimensional case.

The weighted sum scalarization is an objective space method based on the optimization of a weighted sum of the objective functions using non-negative weights. Note that not all efficient solutions can be determined as optimal solutions of the weighted sum scalarization using suitable weights \citep[see, e.\,g.][]{aneja1979bicriteria}. Efficient solutions which can be obtained by using weighted sum scalarization are denoted as \emph{supported efficicent} and their corresponding non-dominated points are located on the boundary of the convex hull of feasible image points.  Extensions of the weighted sum method to the multi-objective case are proposed in \citet{Przybylski2010a}, \citet{oezpeynirci_2010an}, \citet{Boekler2015output}, and \citet{przybylski2019simple}. 

\citet{Ulungu1995two} introduced the so-called two-phase method for bi-objective problems. In the first phase the extreme supported non-dominated points are generated with an algorithm similar to the initial weighted sum approach. In the second phase the remaining non-dominated points are generated by searching in triangles defined by two consecutive extreme supported non-dominated points. In \citet{Przybylski2008two} and \citet{Tuyttens2000performance} problem specific algorithms are suggested for the  second phase, while in \citet{Przybylski2010atwo} a two-phase method for problems with more than two objectives is proposed. 

The augmented weighted Tchebycheff method, first presented in \citet{Steuer1983an}, minimizes the augmented weighted Tchebycheff distance between a predefined reference point and the set of feasible image points. \citet{Daechert2012an} suggested an adaptive choice of the augmentation term for the bi-objective case.

In \citet{Boland2015balanced}, \citet{Boland2015triangle} (for the bi-objective case), \citet{Daechert2014a}, and \citet{Klamroth2015on} (for the tri- respectively multi-objective case) search region splitting methods are proposed. In this class of objective space methods, the search region (based on the already determined non-dominated points) is splitted into so-called search zones on which scalarizations are solved indpendently.

Besides their advantages, objective space methods share a shortcoming: In each iteration a scalarized integer program is solved from scratch. Even though in some objective space methods starting solutions can be transfered from previous iterations, a large number of very similar problems has to be solved. In order to avoid this effort, decision space methods, mainly the branch and bound method, have been increasingly investigated in the recent years.

\citet{Klein1982an} developed one of the first branch and bound algorithms for multi-objective integeger programs with a typical one tree structure. In \citet{Kiziltan1983an} a general branch and bound framework for multi-objective integer programs with binary variables is presented.  \citet{Ulungu1997solving} and \citet{Visee1998two} proposed problem specific branch and bound approaches for bi-objective Knapsack problems, where the latter approach is integrated in a two-phase method.

\citet{Mavrotas1998a} extend the branch and bound approach to multi-objective mixed integer programs. Parts of the algorithm are refined in \citet{Mavrotas2005multi}. In \citet{Vincent2013multiple} this algorithm is improved and it is shown that the original algorithm is not correct because the final dominance test is incomplete. In \citet{Belotti2013a} a branch and bound method is presented that can handle bi-objective mixed integer programs with continious variables in both objective functions.

The branch and bound method proposed in \citet{Sourd2008a} uses a set of points as lower bound instead of just using a single point. Furthermore hyperplanes are used to fathom nodes by dominance. In \citet{Stidsen2014a} this idea is continued. They use hyperplanes as a lower bound set that are generated by solving weighted sum scalarizations. Additionally they present the so-called \emph{Pareto branching} and the \emph{slicing} technique. With Pareto branching it is possible to divide the objective space to possibly ignore parts of it in specific nodes. Slicing partitions the objective space in equally large parts and a respective slice can be fathomed if it is dominated by an already found integer point. In \citet{Stidsen2018a} this algorithm is improved and an approach to parallelize the algorithm is presented. Based on this, the idea Pareto branching is further investigated in \citet{Parragh2019branch} and \citet{Gadegaard2019bi} for the bi-objective case and \citet{Forget2020branch} for the tri-objective case.
A self-contained survey of multi-objective branch and bound approaches is given in \citet{przybylski17multi}. 

In this paper we present a bi-objective branch and bound algorithm that is augmented by scalarization-based information. We make use of optimized single-objective solvers for scalar integer programs and integrate the resulting information into the bi-objective branch and bound by improving lower and upper bounds.
Furthermore, we propose a new adaptive node selection strategy, which relies on objective space information. In our numerical analysis we show the effectiveness of these improvements by comparing them with a generic multi-objective branch and bound algorithm, which we use as our baseline algorithm.

The remainder of the article is organized as follows: In Section~\ref{sec:def}, we introduce notations and definitions for multi-objective optimization. In Section~\ref{sec:BB}, we present a general multi-objective branch and bound framework and its key components. Furthermore, we describe a specific (however standard) multi-objective branch and bound algorithm, which will be used as baseline implementation in our numerical tests. In Section~\ref{sec:hyb}, we present augmentations of the multi-objective branch and bound, that utilize objective space information to improve the node selection as well as the computation of upper and lower bounds.  We provide numerical results in Section~\ref{sec:res} and in Section~\ref{sec:dis}, we outline conclusions and outlooks for further research.

%% file: sec_def.tex
We introduce a general \emph{multi-objective integer linear program} which can be written in the form:
\begin{equation}\label{eq:1} \tag{MOILP}
   \begin{array}{rr@{\extracolsep{1ex}}c@{\extracolsep{1ex}}ll}
      \min  & \multicolumn{3}{@{\extracolsep{0.75ex}}l}{\displaystyle \bigl( z_1(x),\ldots, z_p(x)\bigr)^\top}\\
      \mathrm{ s.t.} & \displaystyle A\,x &\leq&  b  \\
      &x &\geq& 0 \\
      &x &\in& \mathbb{Z}^n .
   \end{array}
\end{equation}
Thereby, \(z(x)\coloneqq (z_1(x),\ldots, z_p(x))^\top = C\cdot x \in\R^p\) (with $p \geq 2$) denotes the objective function vector, with $C \in \R^{p\times n}$ the matrix of objective coefficients.
The set of feasible solutions $X\coloneqq \{x \in\mathbb{Z}^n: A \leq b, x\geq 0 \}$ is a subset of the \emph{decision space} \(\R^n\), while its image $Y\coloneqq \{C\, x\colon x\in X\}$ is a subset of the \emph{objective space} \(\R^p\).

We use the \emph{Pareto concept of optimality} which relies on the componentwise order. 
Let $y^1,y^2 \in \R^p$, then we define the corresponding dominance relations as follows:
\begin{itemize}
\item  $y^1 \leqq y^2$, i.e., $y^1$ \emph{weakly dominates} $y^2$  if $y^1_k \leq y^2_k$ for $k=1,...,p$,
\item  $y^1 < y^2$, i.e., $y^1$ \emph{strictly dominates} $y^2$ if $y^1_k < y^2_k$ for $k=1,...,p$,
\item  $y^1 \leq y^2$, i.e., $y^1$ \emph{dominates} $y^2$ if $y^1 \leq  y^2$ and $y^1 \neq y^2$.
\end{itemize}
A feasible solution $x\in X$ is called \emph{efficient} if there is no other solution $\hat{x}\in X$ dominating it, i.e., $z(\hat{x}) \leq z(x)$. A feasible solution $x\in X$ is called \emph{weakly efficient} if there is no $\hat{x}\in X$ such that $z(\hat{x}) < z(x)$. The set of  efficient solutions is denoted by  $X_E$. By $Y_N= \{ z(x) \in Y\colon x\in X_E \}$ we denote the set of the non-dominated points in the objective space. Moreover, for any set \(Q\subseteq \R^p\) we denote by \(Q_N\) the set of its non-dominated points (i.e., \(q\in Q_N\iff \nexists q'\in Q\colon q'\leq q\)).  For a comprehensive introduction to multi-objective optimization see, e.\,g., \citet{ehrgott05multicriteria}.

In this article we consider a \emph{minimal complete set} as solution of a multi-objective optimization problem. A minimal complete set denotes the set of all non-dominated points  $Y_N$ and one efficient solution for each non-dominated point. See 
\mbox{\citet{Serafini1987some}}
for a comparison of solution concepts in multi-objective optimization.

A standard solution approach in multi-objective optimization is the \emph{weighted sum scalarization} given in \eqref{eq:3}. 

\begin{equation}\label{eq:3}\tag{$\WS_\lambda$}
\begin{split}
 \min \;\; & \WS_\lambda(x)\coloneqq\lambda^\top z(x) =  \sum_{i=1}^p \lambda_i\, z_i(x)\\
 \mathrm{s.t.}\;\; & x \in X
 \end{split}
\end{equation}
Obviously, every optimal solution of the weighted sum scalarization for $\lambda \in \R^p_>\coloneqq\{\lambda\in\R^p\colon \lambda>0 \}$ is efficient for \eqref{eq:1}. However, in general not all efficient solutions are optimal solutions of a corresponding weighted sum problem. An efficient solution \(x'\in X_E\) is called \emph{supported} if there is a weighting vector $\lambda' \in \R^p_>$ such that \(x'\) is optimal for \eqref{eq:3} for \(\lambda=\lambda'\), otherwise \(x'\) is \emph{unsupported}. Note that the non-dominated points corresponding to supported efficient solutions are located on the boundary of the convex hull of $Y$, while the unsupported non-dominated points are located in its (relative) interior.

As already mentioned in the introduction the computation of upper and lower bounds on the non-dominated set is a crucial component of any multi-objective branch and bound algorithm. The tightest componentwise upper and lower bounds of \(Y_N\) are  the \emph{ideal point} $y^I$ and the \emph{Nadir point} $y^N$ given by:
\[
  y^I_k = \min_{y\in Y} y_k  \quad \text{and}\quad y^N_k = \max_{y\in Y_N} y_k \qquad \text{for } k=1,\ldots p. 
\]
Obviously, $y^I \leqq y \leqq y^N$ holds for every $y\in Y_N$, i.e.; \(Y_N\) is contained in the hyperbox spanned by the corner points \(y^I\) and \(y^N\). However, these single point bounds are in general very weak except for the degenerate case of $y^I = y^N$. This motivates to consider \emph{bound sets} instead of bounds consisting of a single point. We will rely on the definition of bound sets proposed in \citet{ehrgott2007bound}. Let $\R^p_\geqq \coloneqq \{y \in \R^p\colon y\geqq 0\}$, then
\begin{itemize}
   \item A \emph{lower bound set}  $\LBS \subset \R^p$ for $Y_N$ is a 
      \begin{itemize}[noitemsep]
         \item $\R^p_\geqq$-closed (i.e., the set \(\LBS +\R^p_\geqq\) is closed),
         \item $\R^p_\geqq$-bounded (i.e., there exists a \(y\in\R^p\) such that \(\LBS\subset y+\R^p_\geqq\)) 
         \item stable set (i.e., $\LBS \subset (\LBS +\R^p_\geqq)_N$),
      \end{itemize}
      such that $Y_N \subset (\LBS +\R^p_\geqq)$.
   \item An \emph{upper bound set} $\UBS \subset \R^p$ for $Y_N$ is a 
      \begin{itemize}[noitemsep]
         \item $\R^p_\geqq$-closed,
         \item $\R^p_\geqq$-bounded,
         \item stable sets,
      \end{itemize}
   such that $Y_N \subset \cl\bigl((\UBS +\R^p_\geqq)^\complement \bigr)$.
\end{itemize}
The upper bound and lower bound that we will define for our branch and bound framework in Section \ref{sec:BB} will suit these definitions. We say a lower bound $\LBS $ is weakly dominated by an upper bound $\UBS$ if for all $l \in \LBS $ there exists an $u\in \UBS $ such that $u \leqq l$.

In the following we restrict ourselves to bi-objective binary linear optimization problems, i.\,e., problems with two linear objective functions and variables \(x \in \{ 0,1 \}^n\):
\begin{equation}\label{eq:2}\tag{BO01LP}
\begin{array}{rr@{\extracolsep{0.75ex}}c@{\extracolsep{0.75ex}}l}
 \min  & \displaystyle z(x) &=& \bigl(z_1(x),z_2(x)\bigr)^\top\\
 \text{ s.t.} & \displaystyle A\,x &\leq& b  \\
              &       x &\in& \{ 0,1 \}^n.
\end{array}
\end{equation}

%% file: sec_BB.tex
In this section we present a generic multi-objective branch and bound framework, which we specify and augment by using scalarization based information in the then following sections.

branch and bound methods follow a ``divide and conquer'' paradigm. A problem that is too hard to be solved directly, is divided into smaller and thus easier subproblems. 
Thereby, subproblems are associated with nodes in a tree data structure according to their descent, i.e., node \(i\) is a descendant node of node \(j\) iff the feasible set of the subprobem associated with node \(i\) is a subset of the feasible set of the subproblem associated with node \(j\). The corresponding subproblems of the child nodes are created by subdividing the feasible set of the corresponding (sub)problem of the parent node. 
Starting with the root node, to which the original optimization problem is associated, the algorithm selects in each iteration one active node and updates its lower bound and upper bound. Then the active node can be fathomed if the corresponding subproblem is either solved or irrelevant for the determination of a minimal complete set. If we cannot prune we subdivide the corresponding problem into new subproblems and create corresponding child nodes (branching). For a more detailed introduction and survey of multi-objective branch and bound algorithms see \citet{przybylski17multi}. A recent survey of single-objective branch and bound frameworks is given e.g.\ in \citet{Morrison2016branch}.
In the following we specify the lower bound, upper bound, branching rule and node selection we use in our framework.

\paragraph{Lower bound:} Lower bound sets are often determined by solving relaxations of the respective subproblem. Like in the single-objective case, the most frequently used relaxations are \emph{linear} and \emph{convex relaxations}. In order to solve the linear relaxation we are using in our framework, we apply \emph{Benson's outer approximation algorithm} \citep{Benson1998an,ehrgott12dual}. The algorithm is initiated with a lower bound, which is improved in every iteration by generating cuts. Due to the outer approximation structure the algorithm can be aborted at any time returning a valid lower bound. Alternatively, linear (or convex) relaxations can be obtained using a dichotomic scheme \citep[see, for example,][]{aneja1979bicriteria,oezpeynirci_2010an,Przybylski2010a}.
In the following we denote a lower bound set $\LBS$ as \emph{convex lower bound set} or \emph{convex lower bound} if $\LBS + \R^2_\geq$ is a convex set. Note that the set \(\LBS\) is thereby not necessarily convex.

\paragraph{Upper bound:} The upper bound set, in the following denoted by \(\UBS\), is stored in the form of a so-called \emph{incumbent list}. Throughout the run of the algorithm, it contains all integer feasible solutions and their corresponding outcome vectors that are not dominated by another feasible solution found so far. 
In every iteration the extreme supported solutions of the computed lower bound sets are checked for integer feasibility. An integer feasible solution $\bar{x}\in X$ is then appended to the incumbent list, if there is no \(x\in\UBS\) dominating \(\bar{x}\), i.\,e., $C(x)\leq C(\bar{x})$. 
If a new solution $\bar{x}$ is added to the incumbent list \(\UBS\) all solutions \(x\in\UBS\) which are dominated by \(\bar{x}\) ($C(\bar{x})\leq C(x)$) are removed from it, that is
\[
\UBS \uplus \{\bar{x}\} \coloneqq 
\begin{cases}
	\UBS & \text{if}\; \exists x\in\UBS\colon C(x)\leq C(\bar{x}) \\
	\{\bar{x}\}\cup \{x\in \UBS \colon C(\bar{x}) \nleq C(x) \} & \text{otherwise}.
\end{cases}
\]
 Note that an update of the incumbent list requires a subsequent update of the list of \emph{local upper bounds}. A detailed description of local upper bounds, their computation and update in an arbitrary number of criteria is given in \citet{Klamroth2015on}. In this framework we start with an empty upper bound set. However, it is also possible to initialize the incumbent list by heuristic methods, or by solving scalarizations like, e.g., in the two-phase method \citep{Ulungu1995two,Visee1998two}.

\paragraph{Node selection:} In every iteration of the algorithm an unexplored node is selected from the tree of subproblems. This node is called \emph{active node}. The order in which the nodes of the tree are considered has a significant impact on the number of created nodes that have to be explored and thus on the computation time. 

Two types of strategies need to be distinguished: \emph{static} strategies and \emph{dynamic} strategies. The  two most common examples of static strategies are the \emph{depth-first} strategy and the \emph{breadth-first} strategy. Most multi-objective branch and bound algorithms in literature follow a depth-first strategy. Thus, we use this strategy for our baseline implementation as well.

In contrast to the single-objective case, dynamic node selection strategies are rarely applied in the multi-objective case. Dynamic node selection strategies are, for example, applied in \citep{Belotti2013a,Stidsen2014a,Jesus2021on}.

\paragraph{Fathoming:} In order to avoid the total enumeration of all  feasible solutions, nodes are fathomed if the respective subproblem is either solved to optimality or does not contain solutions which are necessary to determine a minimal complete set. In particular,  there are three different situations in which a node can be fathomed:

\begin{enumerate}[label=\emph{\roman*})]
   \item \emph{Fathoming by infeasibility:} If the LP-relaxation of a subproblem is infeasible then the corresponding subproblem is infeasible as well, since the feasible set of the subproblem is a subset of the feasible set of its relaxation. 
   \item \emph{Fathoming by optimality:} Similar to the single-objective case we can fathom a node by optimality if the lower bound $\LBS $ is equal to the upper bound $\UBS$. This implies the subproblem is solved to optimality and the associated node must not be subdiveded further. However, this can happen in the multi-objective case only if the lower and upper bound consist of the same single point, namely the ideal point.
   \item \emph{Fathoming by dominance:} A node can be fathomed by dominance if all  feasible solutions of this subproblem are dominated by points in the incumbent list. In order to check dominance for all feasible outcome vectors of a subproblem we compare the lower bound $\LBS $ of the corresponding node to the current upper bound $\UBS $. If for all $l \in \LBS $ there is a point in the incumbent list $u\in \UBS $ with $u \leqq l$ then all feasible points in the subtree are dominated by the current incumbent list. In other words, if there is no local upper bound defined by $\UBS $ above the computed lower bound the node can be fathomed by dominance.
\end{enumerate}

\paragraph{Branching:} As mentioned in the beginning of this section, one of the key aspects of branch and bound is iterative subdivision into smaller subproblems. Thereby subproblems are associated with nodes in a tree, such that the subproblem associated to a child node is obtained by one branching step. Since we consider binary optimization problems (\ref{eq:2}), we can divide a (sub)problem into two new subproblems by fixing a specific variable to $0$ and respectively to $1$ in the other subproblem. This results in a binary branch and bound tree.

The branching rule determines which variable is selected as branching variable in each iteration. Thereby, one distinguishes between static and dynamic strategies. Static strategies determine an order of the variables in advance. In each iteration of the algorithm the next variable in this list is used as branching variable. With dynamic strategies the branching variable is selected by considering information obtained from previous iterations, i.e., from the solution of (linear) relaxations of (sub)problems. 

The basic idea of static strategies for single-objective problems is to sort the variables, beginning with the most promising according to the objective function values (see, e.\,g., \citep{Kellerer2004knapsack}). However, this cannot be easily extended to the multi-objective case due to conflicting objective functions. Nevertheless there are some approaches to extend static strategies to the multi-objective case (see for example \citep{Ulungu1997solving} and \citep{Bazgan2009solving}).

In contrast to most of the published papers which apply static strategies we use a dynamic strategy as proposed in \citet{Belotti2013a}. By solving the linear relaxation of a (sub)problem we obtain the lower bound set $\LBS $. For all extreme points of $\LBS $ we check how often a variable is fractional in the corresponding solutions. As branching variable we choose the one which is most often fractional.

%% file: sec_hyb.tex
In this section, we propose modifications which improve the computational efficiency of bi-objective branch and bound algorithms in two critical aspects. One of the weaknesses of multi-objective branch and bound as compared to its single-objective counterpart is the bounding procedure. While any feasible solution \(\bar{x}\in X\) dominates w.r.t.\ one (linear) objective a half-space in decision space (i.e., \(\{x\in\R^n\colon c^\top x \geq c^\top \bar{x}\}\)), the set of feasible solutions which are dominated by a solution \(\bar{x}\) in \(p\geq 2\) objective functions (\(C\in\R^{p\times n}\)) forms a cone \(\{x\in\R^n\colon C\, x \geqq C \,\bar{x}\}\). The cone of dominated solutions is smaller the more the objective functions are in conflict, leading also to a larger number of efficient solutions. This implies that a significant part of the branch and bound tree has to be enumerated and only a small number of branches can be pruned by dominance. 
Despite of this general problem in multi-objective optimization, this asks for good bounding procedures to avoid the unnecessary evaluation of dominated branches. This however, requires good solutions in the incumbent list as well as tight lower bounds.

In order to achieve this, we suggest a new branching strategy and the hybridization of branch and bound with objective space methods. We determine scalarized subproblems adapted to the state of the branch and bound and solve these to integer optimality.

\subsection{Branching Strategy}\label{sec:branch}
 The branching strategy comprises two subsequent decisions: the choice of the active node and its branching into subproblems, i.\,e.\ the decision on which variable the subproblem is branched. This second step is denoted as branching rule. We discuss these two steps together since the order in which the nodes are considered has a significant impact on the branched variable. Instead of the static depth- or breadth-first we use a dynamic node selection strategy, while we rely on the most fractional rule as branching rule. 

The basic idea of our strategy is quite simple and a natural extension of choosing the largest gap in the single-objective case \citep[see, for example,][]{Dechter1985generalized}. For every created node we compute the \emph{approximate hypervolume gap} between lower and upper bound. We use the definition of hypervolume proposed in \citet{Zitzler1999multiobjective}. In every iteration we choose the node with the largest hypervolume gap as active node (cf.\ \cite{Jesus2021on}). Note that when a node is created during the branching process, the approximated hypervolume gap of the parent node is assigned to it. 
We distinguish two variants of the hypervolume gap: the \emph{total hypervolume gap} and the \emph{local hypervolume gap}. While the total hypervolume gap measures the volume of the search region, i.\,e.\ the volume between lower and upper bound set, the local hypervolume gap approach considers only the volume of the largest search zone, i.\,e.\ the gap between a local upper bound and the lower bound set. For a more detailed definition of search regions and search zones we refer to \citet{Klamroth2015on}. 

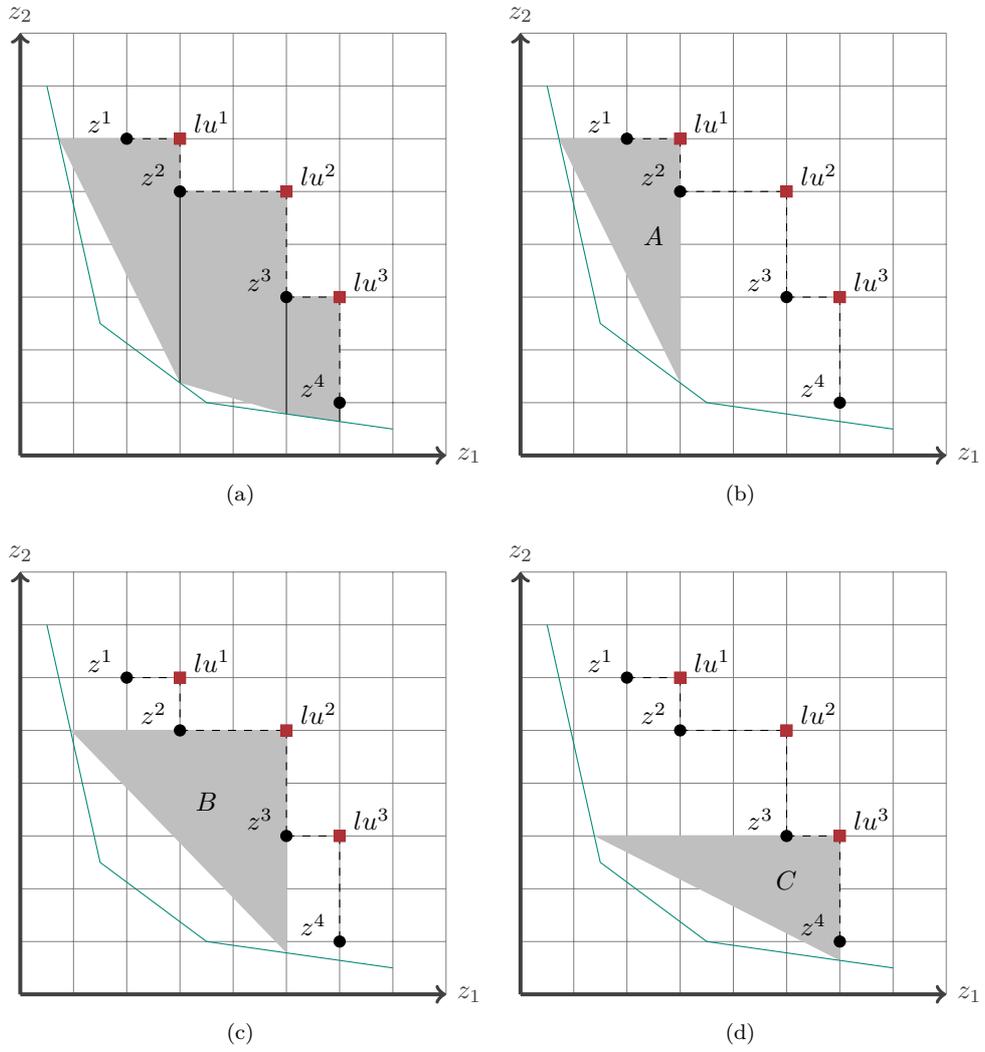
\begin{figure}[htbp!]
   \subfloat[\label{fig1a}]{
		\begin{tikzpicture}[scale=0.7]
		\draw[style=help lines] (0,0) grid (8,8);
		\draw[line width=1.5,->,black!75](0,0) -- (8,0) node[right] {$z_1$};
		\draw[line width=1.5,->,black!75] (0,0) -- (0,8) node[above] {$z_2$};

		\draw[fill,lightgray] (0.7222,6) -- (3,6) -- (3,1.375);
		\draw[fill,lightgray]  (3,1.375) -- (3,5) -- (5,5) -- (5,0.7857);
		\draw[fill,lightgray]  (5,0.7857) -- (5,3) -- (6,3) -- (6,0.642);
		
		\draw (3,5)  -- (3,1.375);
		\draw (5,3) -- (5,0.7857);
		\draw (6,1) -- (6,0.642);

		\draw (3.6,6.7) node[below,black] {$lu^1$};
		\draw (5.6,5.7) node[below,black] {$lu^2$};
		\draw (6.6,3.7) node[below,black] {$lu^3$};		
		\draw (1.5,6.7) node[below,black] {$z^1$};
		\draw (2.5,5.7) node[below,black] {$z^2$};
		\draw (4.5,3.7) node[below,black] {$z^3$};
		\draw (5.5,1.7) node[below,black] {$z^4$};
		
		\filldraw (2,6) circle (3pt);
		\filldraw (6,1) circle (3pt);
		\filldraw (5,3) circle (3pt);
		\filldraw (3,5) circle (3pt);
		
		\draw[PineGreen] (0.5,7) -- (1.5,2.5) -- (3.5,1) -- (7,0.5);
		
		\draw[dashed] (2,6) -- (3,6);
		\draw[dashed] (3,6) -- (3,5);
		\draw[dashed] (3,5) -- (5,5);
		\draw[dashed] (5,5) -- (5,3);
		\draw[dashed] (5,3) -- (6,3);
		\draw[dashed] (6,3) -- (6,1);
		
		\filldraw[Maroon] ([xshift=-3pt,yshift=-3pt]3,6) rectangle ++(6pt,6pt);
		\filldraw[Maroon] ([xshift=-3pt,yshift=-3pt]5,5) rectangle ++(6pt,6pt);
		\filldraw[Maroon] ([xshift=-3pt,yshift=-3pt]6,3) rectangle ++(6pt,6pt);
	\end{tikzpicture}}	
   \subfloat[\label{fig1b}]{
	\begin{tikzpicture}[scale=0.7]
      \draw[style=help lines] (0,0) grid (8,8);
      \draw[line width=1.5,->,black!75](0,0) -- (8,0) node[right] {$z_1$};
      \draw[line width=1.5,->,black!75] (0,0) -- (0,8) node[above] {$z_2$};
      
      \draw[fill,lightgray] (0.7222,6) -- (3,6) -- (3,1.375);
      
      \draw (2.5,4.5) node[below,black] {$A$};
      \draw (3.6,6.7) node[below,black] {$lu^1$};
      \draw (5.6,5.7) node[below,black] {$lu^2$};
      \draw (6.6,3.7) node[below,black] {$lu^3$};
      \draw (1.5,6.7) node[below,black] {$z^1$};
      \draw (2.5,5.7) node[below,black] {$z^2$};
      \draw (4.5,3.7) node[below,black] {$z^3$};
      \draw (5.5,1.7) node[below,black] {$z^4$};
      
      \filldraw (2,6) circle (3pt);
      \filldraw (6,1) circle (3pt);
      \filldraw (5,3) circle (3pt);
      \filldraw (3,5) circle (3pt);
      
      \draw[PineGreen] (0.5,7) -- (1.5,2.5) -- (3.5,1) -- (7,0.5);
            
      \draw[dashed] (2,6) -- (3,6);
      \draw[dashed] (3,6) -- (3,5);
      \draw[dashed] (3,5) -- (5,5);
      \draw[dashed] (5,5) -- (5,3);
      \draw[dashed] (5,3) -- (6,3);
      \draw[dashed] (6,3) -- (6,1);
      
      \filldraw[Maroon] ([xshift=-3pt,yshift=-3pt]3,6) rectangle ++(6pt,6pt);
      \filldraw[Maroon] ([xshift=-3pt,yshift=-3pt]5,5) rectangle ++(6pt,6pt);
      \filldraw[Maroon] ([xshift=-3pt,yshift=-3pt]6,3) rectangle ++(6pt,6pt);
   \end{tikzpicture}}
   
   \subfloat[\label{fig1c}]{
	\begin{tikzpicture}[scale=0.7]
		\draw[style=help lines] (0,0) grid (8,8);
		\draw[line width=1.5,->,black!75](0,0) -- (8,0) node[right] {$z_1$};
		\draw[line width=1.5,->,black!75] (0,0) -- (0,8) node[above] {$z_2$};
		
      \draw[fill,lightgray] (0.9444,5 ) -- (5,5) -- (5,0.7857);

		\draw (3.5,4) node[below,black] {$B$};
			\draw (3.6,6.7) node[below,black] {$lu^1$};
		\draw (5.6,5.7) node[below,black] {$lu^2$};
		\draw (6.6,3.7) node[below,black] {$lu^3$};
		\draw (1.5,6.7) node[below,black] {$z^1$};
		\draw (2.5,5.7) node[below,black] {$z^2$};
		\draw (4.5,3.7) node[below,black] {$z^3$};
		\draw (5.5,1.7) node[below,black] {$z^4$};
		
		\filldraw (2,6) circle (3pt);
		\filldraw (6,1) circle (3pt);
		\filldraw (5,3) circle (3pt);
		\filldraw (3,5) circle (3pt);

		\draw[PineGreen] (0.5,7) -- (1.5,2.5) -- (3.5,1) -- (7,0.5);

		\draw[dashed] (2,6) -- (3,6);
		\draw[dashed] (3,6) -- (3,5);
		\draw[dashed] (3,5) -- (5,5);
		\draw[dashed] (5,5) -- (5,3);
		\draw[dashed] (5,3) -- (6,3);
		\draw[dashed] (6,3) -- (6,1);
		
		\filldraw[Maroon] ([xshift=-3pt,yshift=-3pt]3,6) rectangle ++(6pt,6pt);
		\filldraw[Maroon] ([xshift=-3pt,yshift=-3pt]5,5) rectangle ++(6pt,6pt);
		\filldraw[Maroon] ([xshift=-3pt,yshift=-3pt]6,3) rectangle ++(6pt,6pt);
	\end{tikzpicture}}
   \subfloat[\label{fig1d}]{
   \begin{tikzpicture}[scale=0.7]
      \draw[style=help lines] (0,0) grid (8,8);
      \draw[line width=1.5,->,black!75](0,0) -- (8,0) node[right] {$z_1$};
      \draw[line width=1.5,->,black!75] (0,0) -- (0,8) node[above] {$z_2$};
      
      \draw[fill,lightgray] (6,0.642) -- (6,3) -- (1.388, 3);
      
      \draw (5,2.5) node[below,black] {$C$};
     	\draw (3.6,6.7) node[below,black] {$lu^1$};
      \draw (5.6,5.7) node[below,black] {$lu^2$};
      \draw (6.6,3.7) node[below,black] {$lu^3$};
      \draw (1.5,6.7) node[below,black] {$z^1$};
      \draw (2.5,5.7) node[below,black] {$z^2$};
      \draw (4.5,3.7) node[below,black] {$z^3$};
      \draw (5.5,1.7) node[below,black] {$z^4$};

      \filldraw (2,6) circle (3pt);
      \filldraw (6,1) circle (3pt);
      \filldraw (5,3) circle (3pt);
      \filldraw (3,5) circle (3pt);
      
		\draw[PineGreen] (0.5,7) -- (1.5,2.5) -- (3.5,1) -- (7,0.5);
      
      \draw[dashed] (2,6) -- (3,6);
      \draw[dashed] (3,6) -- (3,5);
      \draw[dashed] (3,5) -- (5,5);
      \draw[dashed] (5,5) -- (5,3);
      \draw[dashed] (5,3) -- (6,3);
      \draw[dashed] (6,3) -- (6,1);
      
      \filldraw[Maroon] ([xshift=-3pt,yshift=-3pt]3,6) rectangle ++(6pt,6pt);
      \filldraw[Maroon] ([xshift=-3pt,yshift=-3pt]5,5) rectangle ++(6pt,6pt);
      \filldraw[Maroon] ([xshift=-3pt,yshift=-3pt]6,3) rectangle ++(6pt,6pt);
   \end{tikzpicture}}
	\caption{Example of computation of the two different approximated hypervolume gap approaches.}
	\label{fig:hyperv}
\end{figure}

Figure \ref{fig:hyperv} illustrates the two different approaches. Here, $z^1,\ldots, z^4 \in K \subset\UBS$ are points of the incumbent list and $lu^1,\ldots, lu^3$ are their corresponding local upper bounds, where $K$ is a subset of the incumbent list containing just the points above the lower bound of node $\bar{n}$.  The green line represents the lower bound. Figure \ref{fig1a} shows how to measure the total hypervolume gap of a node $\bar{n}$, in the following denoted by $thg(\bar{n})$.  For this approach we consider the approximated search region of the corresponding node. Since there is a natural order in the bi-objective case, it is possible to consider the approximated search zone of the first local upper bound,  i.e. the local upper bound with the smallest $z_1$-value. Therefore we define the two spanning points, which, together with the corresponding local upper bound, define a triangle. The spanning points of a local upper bound $lu$ are defined by \(sp^i (lu) \coloneqq \{ l \in \LBS \colon l_{3-i} = lu_{3-i} \}, i = 1,2\).  So, the approximate hypervolume gap of $lu$ is given by
\[ hg(lu) \coloneqq \frac{1}{2} \,  \bigl\vert sp^1(lu)_1 - lu_1\bigr\vert \cdot \bigl\vert sp^2(lu)_2 - lu_2\bigr\vert.  \]
For the remaining local upper bounds we compute the hypervolume of slices as shown in the illustration. The hypervolume of the slice of $lu^i,i=1,\dots,  \vert K \vert-1$ is defined as
\[sl(lu^i)\coloneqq \frac{\bigl\vert z^i_2 - sp^2(lu^{i-1})_2\bigr\vert+ \bigl\vert lu^i_2 - sp^2(lu^{i})_2\bigr\vert}{2} \cdot  \bigl\vert z^i_1 - lu^i_1\bigr\vert . \] 
So, the total (approximated) hypervolume gap, which is assigned to node $\bar{n}$, is given by
\[thg(\bar{n}) \coloneqq hg(lu^1) + sl(lu^2)  + \ldots + sl(lu^{\vert K \vert-1}). \]

The Figures \ref{fig1b}, \ref{fig1c} and \ref{fig1d} show the computation of the local hypervolume gap. The local hypervolume gap of a node $\bar{n}$ is considered as the largest approximated hypervolume gap of a local upper bound corresponding to points in $K$. Therefore, the the local hypervolume gap, which is assigned to node $\bar{n}$, is defined by
\[
   lhg(\bar{n}) \coloneqq \max_{i=1,\dots, \vert K\vert -1}hg(lu^i).
\]
In the given example, $B$ is the largest approximated hypervolume and therefore is assigned to node $\bar{n}$.

Note that in our presented algorithms in Section \ref{subsec:AugBB} the local upper bound is initialized with the point $(\infty,\infty)^\top$. Therefore, it is possible to apply the new branching strategies immediately at the beginning of the algorithm. Obviously, this approximation may neglect significantly large parts of the search regions and search zones. However, the idea of the approximated hypervolume gap eases computation and saves time. The efficiency of these new dynamic branching strategies is shown in Section \ref{sec:res}.

\subsection{Augmenting Branch and Bound with IP Scalarizations} \label{subsec:AugBB}
In this subsection, we introduce a method to incorporate scalarizations into branch and bound. We build a hybrid branch and bound algorithm combining the partial enumeration of decision space with objective space information by solving scalarizations to integer optimality.

An integer optimal solution \(\bar{x}\) of a scalarization can be used to update upper and lower bound. Obviously, the corresponding image point \(z(\bar{x})\) can be added to the incumbent list. Moreover,  a scalarizing function and its optimal solution \(\bar{x}\) define a level set, which can be included in the lower bound set for all descendant nodes. In order to utilize these improved lower bounds in all nodes we solve the IP scalarizations in the root node.

\subsubsection{Using Weighted Sum Scalarization}
During the run of the branch and bound algorithm, a strategy triggers the IP solution of weighted sum scalarizations in the root node. Thus, we solve problem \eqref{eq:3} for for adaptively chosen values of $\lambda \in \mathbb{R}^2_>$. Although we solve the IP scalarization in the root node the parameter $\lambda$ is gained from the currently active node.
Thereby, \(\lambda\) is determined by the largest approximated local hypervolume gap in the active node. This gap is spanned by two points in the incumbent list together with their local upper bound. Note that these points spanning the largest gap are already determined if the local hypervolume gap branching strategy is applied. The corresponding value of $\lambda$ is determined by computing the normal to the hyperplane that is defined by those two points.
Once $\lambda$ is obtained, we can  solve problem \eqref{eq:3} with a single-objective integer linear programming solver. Let \(\bar{x}^\lambda\) be the optimal solution of the weighted sum scalarization with weighting vector \(\lambda\), then \(z(\bar{x}^\lambda)\)  is a supported non-dominated point of \eqref{eq:2}. 
Thus, we can add this point to the incumbent list (if it was not found in previous iterations) and filter the resulting list for non-dominance. 
Moreover, the solution of integer scalarizations can also be used to tighten the lower bound set, since the level set \(\{z\in\R^2\colon \lambda^\top z=\WS_\lambda(\bar{x}^\lambda)\}\) provides the valid inequality \(\lambda^\top z(x) \geq \WS_\lambda(\bar{x}^\lambda) \) for all \(x\in X\).

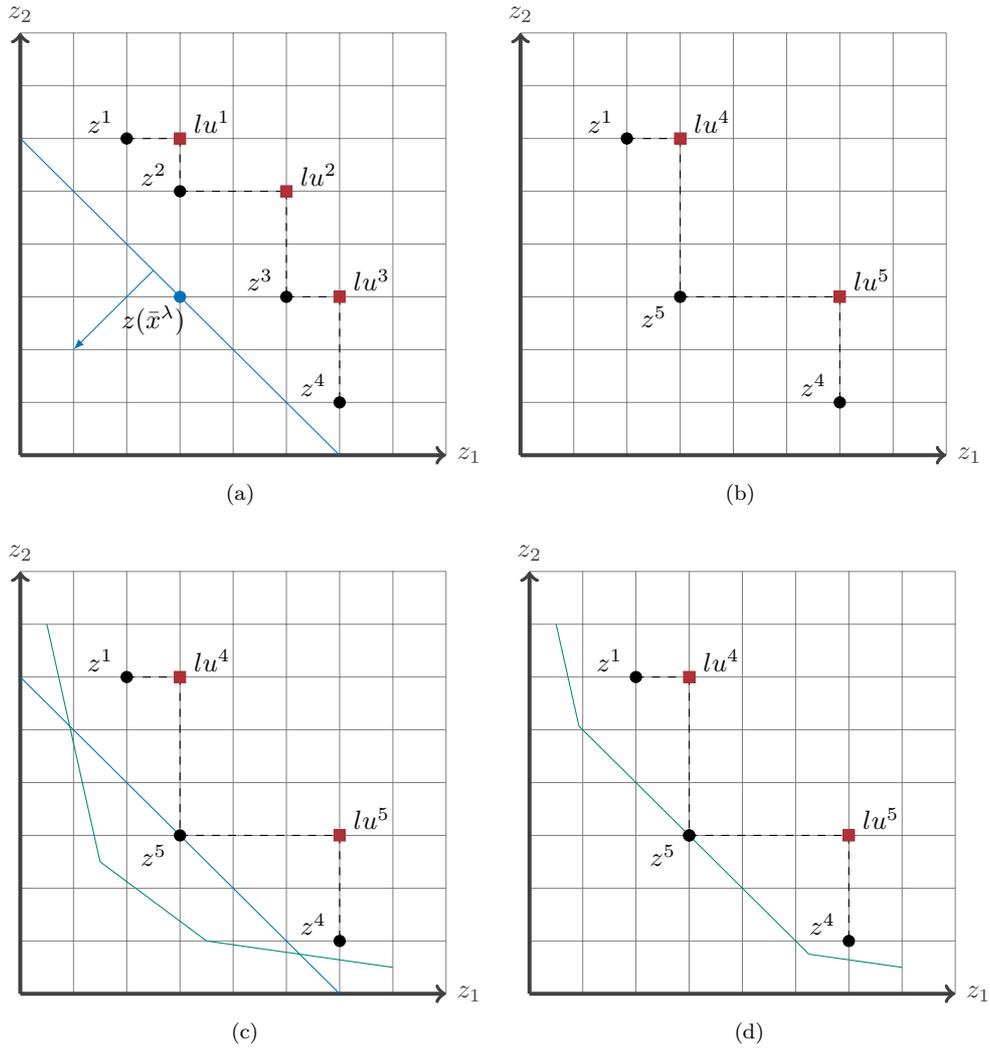
\begin{figure}[h!]
   \subfloat[\label{fig2a}]{
   \begin{tikzpicture}[scale=0.7]

   \draw[style=help lines] (0,0) grid (8,8);
   \draw[line width=1.5,->,black!75](0,0) -- (8,0) node[right] {$z_1$};
   \draw[line width=1.5,->,black!75] (0,0) -- (0,8) node[above] {$z_2$};

   \draw (3.6,6.7) node[below,black] {$lu^1$};
   \draw (5.6,5.7) node[below,black] {$lu^2$};
   \draw (6.6,3.7) node[below,black] {$lu^3$};
   \draw (1.5,6.7) node[below,black] {$z^1$};
   \draw (2.5,5.7) node[below,black] {$z^2$};
   \draw (4.5,3.7) node[below,black] {$z^3$};
   \draw (5.5,1.7) node[below,black] {$z^4$};
   \draw (2.5,3) node[below,black] {$z(\bar{x}^\lambda)$};

   \draw[RoyalBlue] (0,6) -- (6,0);
   \draw[-latex,RoyalBlue] (2.5,3.5) -- (1,2);

   \filldraw (2,6) circle (3pt);
   \filldraw (6,1) circle (3pt);
   \filldraw (5,3) circle (3pt);
   \filldraw (3,5) circle (3pt);
   \filldraw[RoyalBlue] (3,3) circle (3pt);

   \draw[dashed] (2,6) -- (3,6);
   \draw[dashed] (3,6) -- (3,5);
   \draw[dashed] (3,5) -- (5,5);
   \draw[dashed] (5,5) -- (5,3);
   \draw[dashed] (5,3) -- (6,3);
   \draw[dashed] (6,3) -- (6,1);
   \filldraw[Maroon] ([xshift=-3pt,yshift=-3pt]3,6) rectangle ++(6pt,6pt);
   \filldraw[Maroon] ([xshift=-3pt,yshift=-3pt]5,5) rectangle ++(6pt,6pt);
   \filldraw[Maroon] ([xshift=-3pt,yshift=-3pt]6,3) rectangle ++(6pt,6pt);

   \end{tikzpicture}}
   \subfloat[\label{fig2b}]{
   \begin{tikzpicture}[scale=0.7]
   \draw[style=help lines] (0,0) grid (8,8);
   \draw[line width=1.5,->,black!75](0,0) -- (8,0) node[right] {$z_1$};
   \draw[line width=1.5,->,black!75] (0,0) -- (0,8) node[above] {$z_2$};

   \draw (3.6,6.7) node[below,black] {$lu^4$};
   \draw (6.6,3.7) node[below,black] {$lu^5$};
   \draw (1.5,6.7) node[below,black] {$z^1$};
   \draw (2.5,3) node[below,black] {$z^5$};
   \draw (5.5,1.7) node[below,black] {$z^4$};

   \filldraw (2,6) circle (3pt);
   \filldraw (6,1) circle (3pt);

   \filldraw (3,3) circle (3pt);

   \draw[dashed] (2,6) -- (3,6);
   \draw[dashed] (3,6) -- (3,3);
   \draw[dashed] (3,3) -- (6,3);
   \draw[dashed] (6,3) -- (6,1);
   \filldraw[Maroon] ([xshift=-3pt,yshift=-3pt]3,6) rectangle ++(6pt,6pt);
   \filldraw[Maroon] ([xshift=-3pt,yshift=-3pt]6,3) rectangle ++(6pt,6pt);

   \end{tikzpicture}}

   \subfloat[\label{fig2c}]{
   \begin{tikzpicture}[scale=0.7]
   \draw[style=help lines] (0,0) grid (8,8);
   \draw[line width=1.5,->,black!75](0,0) -- (8,0) node[right] {$z_1$};
   \draw[line width=1.5,->,black!75] (0,0) -- (0,8) node[above] {$z_2$};

   \draw (3.6,6.7) node[below,black] {$lu^4$};
   \draw (6.6,3.7) node[below,black] {$lu^5$};
   \draw (1.5,6.7) node[below,black] {$z^1$};
   \draw (2.5,3) node[below,black] {$z^5$};
   \draw (5.5,1.7) node[below,black] {$z^4$};

   \draw[RoyalBlue] (0,6) -- (6,0);

   \filldraw (2,6) circle (3pt);
   \filldraw (6,1) circle (3pt);

   \filldraw (3,3) circle (3pt);

   \draw[PineGreen] (0.5,7) -- (1.5,2.5);
   \draw[PineGreen] (3.5,1) -- (1.5,2.5);
   \draw[PineGreen] (3.5,1) -- (7,0.5);

   \draw[dashed] (2,6) -- (3,6);
   \draw[dashed] (3,6) -- (3,3);
   \draw[dashed] (3,3) -- (6,3);
   \draw[dashed] (6,3) -- (6,1);
   \filldraw[Maroon] ([xshift=-3pt,yshift=-3pt]3,6) rectangle ++(6pt,6pt);
   \filldraw[Maroon] ([xshift=-3pt,yshift=-3pt]6,3) rectangle ++(6pt,6pt);

   \end{tikzpicture}
   }
   \subfloat[\label{fig2d}]{
   \begin{tikzpicture}[scale=0.7]
   \draw[style=help lines] (0,0) grid (8,8);
   \draw[line width=1.5,->,black!75](0,0) -- (8,0) node[right] {$z_1$};
   \draw[line width=1.5,->,black!75] (0,0) -- (0,8) node[above] {$z_2$};

   \draw (3.6,6.7) node[below,black] {$lu^4$};
   \draw (6.6,3.7) node[below,black] {$lu^5$};
   \draw (1.5,6.7) node[below,black] {$z^1$};
   \draw (2.5,3) node[below,black] {$z^5$};
   \draw (5.5,1.7) node[below,black] {$z^4$};

   \filldraw (2,6) circle (3pt);
   \filldraw (6,1) circle (3pt);

   \filldraw (3,3) circle (3pt);

   \draw[dashed] (2,6) -- (3,6);
   \draw[dashed] (3,6) -- (3,3);
   \draw[dashed] (3,3) -- (6,3);
   \draw[dashed] (6,3) -- (6,1);
   \filldraw[Maroon] ([xshift=-3pt,yshift=-3pt]3,6) rectangle ++(6pt,6pt);
   \filldraw[Maroon] ([xshift=-3pt,yshift=-3pt]6,3) rectangle ++(6pt,6pt);

   \draw[PineGreen] (0.5,7) -- (0.9286,5.0714);
   \draw[PineGreen] (5.25,0.75) -- (0.9286,5.0714);
   \draw[PineGreen] (5.25,0.75) -- (7,0.5);

   \filldraw (3,3) circle (3pt);
   \end{tikzpicture}}
\caption{Example of updating the lower and upper bound  with the usage of the weighted sum scalarization.}
\label{fig:boundupdate}
\end{figure}

Figure \ref{fig:boundupdate} illustrates the update of the lower and upper bound set. In Figure~\ref{fig2a}, $z^1,\dots,z^4$  indicate points that are currently in the incumbent list $\UBS$¸ and $lu^1,\dots,lu^3$ are the corresponding local upper bounds. The point \(z(\bar{x}^\lambda)\) is obtained by solving a  weighted sum scalarization \eqref{eq:3} to integer optimality. Since the new point is not contained in the incumbent list so far, we can update the upper bound as it is shown in Figure \ref{fig2b}. The new incumbent list then reads as \(\UBS\coloneqq \{z(\bar{x}^\lambda)\}\cup\{z\in\UBS\colon z(\bar{x}^\lambda) \nleqq z\}\). Moreover, the lower bound set $\LBS$ can be updated by integrating the blue hyperplane into the lower bound set, i.\,e.\ 
\(\LBS\coloneqq \{z\in \LBS +\R^2_\geqq \colon \lambda^\top z \geq \WS(\bar{x}^\lambda) \}_N\) as it is shown in Figure \ref{fig2c} and \ref{fig2d}. 
In this situation, both ---the lower and upper bound--- are updated, which is not the case in general.

The example illustrates the benefits of hybridizing multi-objective branch and bound with IP scalarizations. Due to weak bounding, nodes may not be fathomed by dominance even if they do not contain additional non-dominated points. The tighter upper bound increases the probability of fathoming a node by dominance in later iterations of the algorithm. Also, the lower bound might be improved. Since we are solving an IP scalarization in the root node, the obtained optimal level set is a valid inequality for all subproblems. 
We combine our new branching strategy and the augmentation with IP scalarizations to our first hybrid branch and bound approach.

\paragraph{Hybrid Branch and Bound Algorithm using Weighted Sum Scalarization}
\begin{itemize}[noitemsep]
\item \emph{Lower bound:} linear relaxation
\item \emph{Upper bound:} incumbent list
\item \emph{Node selection:} node with the largest total/local  hypervolume gap
\item \emph{Branching rule:} most fractional
\item Adaptively solve weighted sum scalarizations in the root node to integer optimality to improve lower and upper bounds by objective space information
\end{itemize}

Instead of using a static depth-first strategy (as in the general branch and bound framework in Section \ref{sec:BB}) we apply the dynamic strategy based on the hypervolume gap (c.f.~Section~\ref{sec:branch}). Even though the extreme points of the lower bound sets might be updated by the weighted sum scalarization, the branching variable is selected based on the original lower bounds. This is due to the fact that the preimages of such intersection points of IP scalarizations and the lower bound set are in general not available.
Note that the weighted sum IP scalarizations are included \emph{adaptively} into the branch and bound. The description of their algorithmic control, however, is postponed to Section~\ref{subsec:contr}. 

In order to conclude the description of the proposed hybrid branch and bound algorithm using weighted sum scalarizations, we want to briefly discuss its advantages and shortcomings. Firstly, it is easy to determine the scalarization parameter $\lambda$ and to integrate the hyperplane into the lower bound set. Its advantage, however, is that the lower bound remains convex. Therefore, the check for fathoming by dominance remains intuitive. Unfortunately, the weighted sum scalarization can only find supported efficient solutions and the lower bound cannot be improved beyond the convex hull of $Y_N$. This motivates us to consider the augmented weighted Tchebycheff scalarization, a scalarization approach which can determine also unsupported efficient solutions.

\subsubsection{Using Augmented Weighted Tchebycheff Scalarization}
We start by defining the weighted Tchebycheff norm: Let $w_i > 0,\; i=1,\ldots,p$ be positive weights with $\sum_{i=1}^p w_i = 1$. Then the \emph{weighted Tchebycheff norm} of a vector \(z\in\R^p\) is defined by
\begin{equation}\label{eq:4}
\Vert z\Vert_\infty^w := \max_{i=1,\ldots,p} \bigl\{ w_i\,\vert z_i\vert \bigr\}\,.
\end{equation}
So, the weighted Tchebycheff scalarization of a multi-objective optimization problem \eqref{eq:1} with respect to a given reference point \(s\in\R^p\) can be written as:
\begin{equation}\label{eq:6}
 \min \bigl\{\Vert z(x)-s \Vert^w_\infty \colon x \in X \bigr\}.
\end{equation}
If the reference point is chosen such that \(s<z(x)\) for all \(x\in X\), every efficient solution can be determined as optimal solution of 
the weighted Tchebycheff scalarization \eqref{eq:6} by variation of \(w\in\R^p_+\) \citep[see, e.g.,][]{miettinen98nonlinear}.
Nevertheless, optimal solutions of the weighted Tchebycheff scalarization correspond in general only to weakly efficient solutions of the multi-objective problem \citep{Steuer1983an,miettinen98nonlinear}. This shortcoming is compensated by an additive augmentation term in the \emph{augmented weighted Tchebycheff norm} 
\begin{equation}
\Vert z\Vert_\tau^w := \Vert z\Vert_\infty^w + \tau\,\Vert z\Vert _1 \,,
\end{equation}
where \(\Vert z\Vert_1 = \vert z_1\vert+\ldots+\vert z_p\vert\) denotes the \(L_1\)-norm, \(w_i \geq 0\), \(i=1,\ldots, p\), \(\sum_{i=1}^p w_i = 1\) and \(\tau>0\).
\citet{Steuer1983an} proposed the \emph{augmented weighted Tchebycheff scalarization} given in (\ref{eq:7}).

\begin{equation}\label{eq:7}\tag{$\AWT_\tau^w$}
	\begin{split}
		\min \;\; & \AWT_\tau^w(x)\coloneqq \Vert z(x)-s\Vert_\tau^w \\
		\mathrm{s.t.}\;\; & x \in X
	\end{split}
\end{equation}
Thereby, the augmentation term makes the augmented weighted Tchebycheff norm a strongly monotone norm and thus the objective function of \eqref{eq:7} a strongly increasing achievement scalarizing function \citep{miettinen98nonlinear}. Consequently, every optimal solution of \eqref{eq:7} is efficient for \eqref{eq:1}.

Note that an appropriate choice of the parameter $\tau$ is difficult in general. On the one hand, too small values of $\tau$ may lead to numerical difficulties. On the other hand, non-supported efficient solutions might be suboptimal for \eqref{eq:7} if the  value of \(\tau\) is too large. However, for bi-objective integer programming \cite{Daechert2012an} propose an adaptive method to determine an optimal value of \(\tau\). We use the proposed parameters $w_1,w_2$ and $\tau$ for our method.

As a reference point $s$ we use the local ideal point of two adjacent non-dominated points. Since the augmented weighted Tchebycheff scalarization can only determine non-dominated points (and the corresponding efficient solutions) which are (strictly) dominated by the reference point, we obtain a non-dominated point in this box.

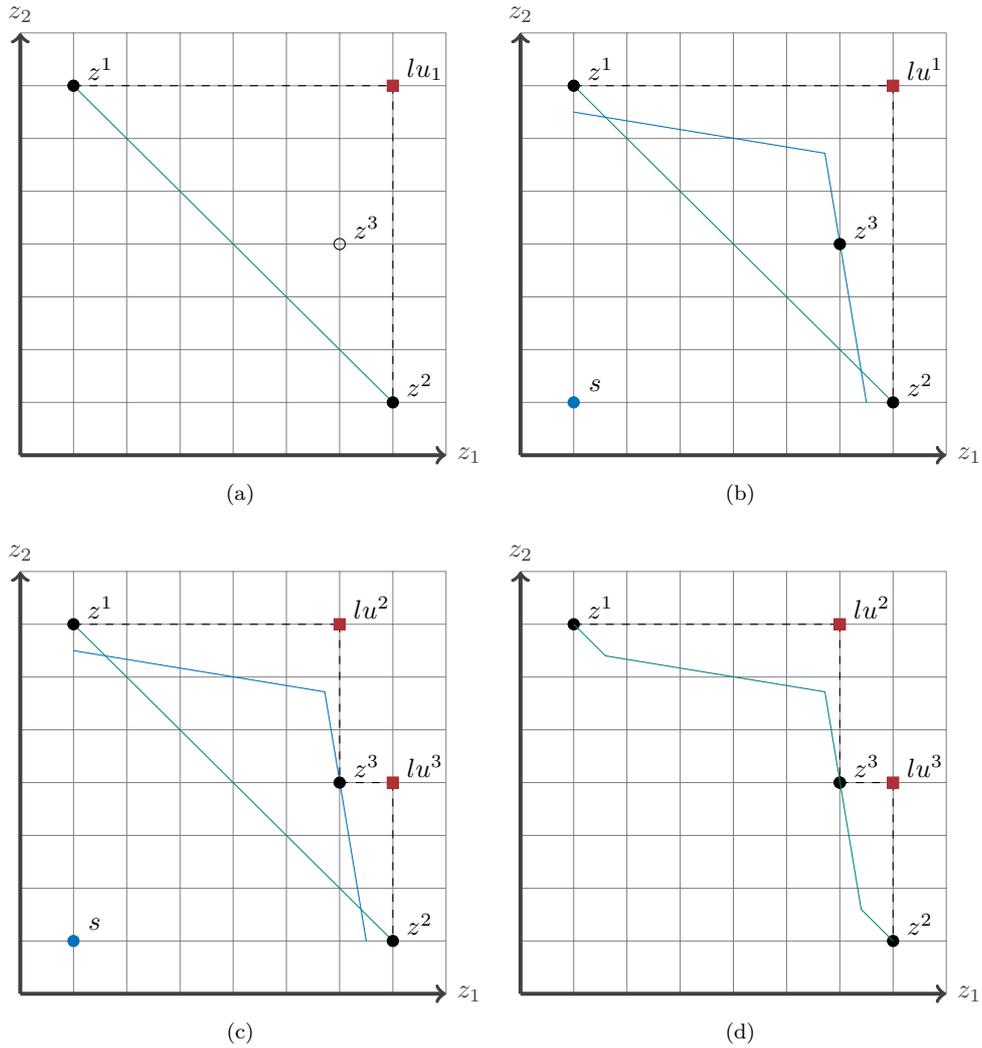
\begin{figure}[h!]\label{fig:awt}

\subfloat[\label{fig3a}]{
 \begin{tikzpicture}[scale=0.7]

\draw[style=help lines] (0,0) grid (8,8);
\draw[line width=1.5,->,black!75](0,0) -- (8,0) node[right] {$z_1$};
\draw[line width=1.5,->,black!75] (0,0) -- (0,8) node[above] {$z_2$};

\draw[PineGreen] (1,7) -- (7,1);

\draw[dashed] (1,7) -- (7,7);
\draw[dashed] (7,7) -- (7,1);

\filldraw[Maroon] ([xshift=-3pt,yshift=-3pt]7,7) rectangle ++(6pt,6pt);

\filldraw (1,7) circle (3pt);
\filldraw (7,1) circle (3pt);
\draw (6,4) circle (3pt);

\draw (1.5,7.7) node[below,black] {$z^1$};
\draw (7.5,1.7) node[below,black] {$z^2$};
\draw (6.5,4.7) node[below,black] {$z^3$};
\draw (7.6,7.7) node[below,black] {$lu_1$};

\end{tikzpicture}}\subfloat[\label{fig3b}]{
 \begin{tikzpicture}[scale=0.7]

\draw[style=help lines] (0,0) grid (8,8);
\draw[line width=1.5,->,black!75](0,0) -- (8,0) node[right] {$z_1$};
\draw[line width=1.5,->,black!75] (0,0) -- (0,8) node[above] {$z_2$};

\draw (1.5,7.7) node[below,black] {$z^1$};
\draw (7.5,1.7) node[below,black] {$z^2$};
\draw (6.5,4.7) node[below,black] {$z^3$};
\draw (7.6,7.7) node[below,black] {$lu^1$};
\draw (1.4,1.6) node[below,black] {$s$};

\draw[PineGreen] (1,7) -- (7,1);

\draw[dashed] (1,7) -- (7,7);
\draw[dashed] (7,7) -- (7,1);

\filldraw[Maroon] ([xshift=-3pt,yshift=-3pt]7,7) rectangle ++(6pt,6pt);

\filldraw (1,7) circle (3pt);
\filldraw (7,1) circle (3pt);
\filldraw[color=RoyalBlue] (1,1) circle (3pt);

\draw[RoyalBlue] (6.5,1) -- (5.72,5.72) -- (1,6.5) ;

\filldraw (6,4) circle (3pt);

\end{tikzpicture}}

\subfloat[\label{fig3c}]{
 \begin{tikzpicture}[scale=0.7]

\draw[style=help lines] (0,0) grid (8,8);
\draw[line width=1.5,->,black!75](0,0) -- (8,0) node[right] {$z_1$};
\draw[line width=1.5,->,black!75] (0,0) -- (0,8) node[above] {$z_2$};

\draw (1.5,7.7) node[below,black] {$z^1$};
\draw (7.5,1.7) node[below,black] {$z^2$};
\draw (6.5,4.7) node[below,black] {$z^3$};
\draw (6.6,7.7) node[below,black] {$lu^2$};
\draw (7.6,4.7) node[below,black] {$lu^3$};
\draw (1.4,1.6) node[below,black] {$s$};

\draw[PineGreen] (1,7) -- (7,1);

\filldraw (1,7) circle (3pt);
\filldraw (7,1) circle (3pt);

\draw[RoyalBlue] (6.5,1) -- (5.72,5.72) -- (1,6.5) ;

\filldraw (6,4) circle (3pt);
\filldraw[color=RoyalBlue] (1,1) circle (3pt);

\draw[dashed] (1,7) -- (6,7);
\draw[dashed] (6,7) -- (6,4);
\draw[dashed] (6,4) -- (7,4);
\draw[dashed] (7,4) -- (7,1);

\filldraw[Maroon] ([xshift=-3pt,yshift=-3pt]7,4) rectangle ++(6pt,6pt);
\filldraw[Maroon] ([xshift=-3pt,yshift=-3pt]6,7) rectangle ++(6pt,6pt);

\end{tikzpicture}}
\subfloat[\label{fig3d}]{
 \begin{tikzpicture}[scale=0.7]

\draw[style=help lines] (0,0) grid (8,8);
\draw[line width=1.5,->,black!75](0,0) -- (8,0) node[right] {$z_1$};
\draw[line width=1.5,->,black!75] (0,0) -- (0,8) node[above] {$z_2$};

\filldraw (1,7) circle (3pt);
\filldraw (7,1) circle (3pt);
\draw (6,4) circle (3pt);

\draw (1.5,7.7) node[below,black] {$z^1$};
\draw (7.5,1.7) node[below,black] {$z^2$};
\draw (6.5,4.7) node[below,black] {$z^3$};
\draw (6.6,7.7) node[below,black] {$lu^2$};
\draw (7.6,4.7) node[below,black] {$lu^3$};

\filldraw (6,4) circle (3pt);

\draw[dashed] (1,7) -- (6,7);
\draw[dashed] (6,7) -- (6,4);
\draw[dashed] (6,4) -- (7,4);
\draw[dashed] (7,4) -- (7,1);

\filldraw[Maroon] ([xshift=-3pt,yshift=-3pt]7,4) rectangle ++(6pt,6pt);
\filldraw[Maroon] ([xshift=-3pt,yshift=-3pt]6,7) rectangle ++(6pt,6pt);

\draw[PineGreen] (1,7) -- (1.599,6.401) -- (5.72,5.72) -- (6.401,1.599) -- (7,1);
\end{tikzpicture}}

\caption{Example of updating the lower and upper bound with the usage of the augmented weighted Tchebycheff scalarization.}
\label{fig:awtch}
\end{figure}

The goal to improve the lower bound set beyond the convex hull of non-dominated points is the motivation to solve augmented weighted Tchebycheff scalarizations to integer optimality. Figure \ref{fig:awtch}  shows an example how such an update of the bounds could look like.  Here, $z^1$ and $z^2$ are two known non-dominated points (obtained with the weighted sum IP scalarization). Point $z^3$ is a non-supported non-dominated point that has not been found yet in Figure \ref{fig3a}.  By using the local ideal point of $z^1$ and $z^2$ as the reference point $s$, Figure \ref{fig3b} illustrates how the non-dominated point $z^3$ is found by applying the augmented weighted Tchebycheff scalarization. In Figure \ref{fig3c} and \ref{fig3d} the resulting improvements of the lower and upper bound are shown. Obviously the lower bound is improved beyond the convex hull of $Y_N$. We now define our second hybrid branch and bound approach:

\paragraph{Hybrid Branch and Bound Algorithm using Augmented Weighted Tchebycheff Scalarization}
\begin{itemize}[noitemsep]
\item \emph{Lower bound:} linear relaxation
\item \emph{Upper bound:} incumbent list
\item \emph{Node selection:} node with the biggest total/local hypervolume gap
\item \emph{Branching rule:} most fractional
\item Adaptively solve weighted sum and augmented weighted Tchebycheff scalarizations in the root node to integer optimality to improve lower and upper bounds by objective space information
\end{itemize}

In addition to the weighted sum scalarization, we use the augmented weighted Tchebycheff scalarization. Since two adjacent non-dominated points are required as input of the augmented weighted Tchebycheff scalarization, we cannot rely on points in the incumbent list, which are only non-dominated so far. In fact, we apply augmented weighted Tchebycheff IP scalarizations only to boxes spanned by points obtained as optimal solutions of the weighted sum scalarization. Thus, we do not rely on parameters from the currently active node, but solve the augmented weighted Tchebycheff scalarization in the largest area defined by two adjacent known non-dominated points. 

When using augmented weighted Tchebycheff IP scalarizations, the lower bound can become tighter than the convex hull of the set of non-dominated points, which reduces the area where new non-dominated points can be found. Additionally, we can find non-supported non-dominated points in early stages of the algorithm. This improves the upper bound in the beginning resulting in a higher chance of fathoming a node by dominance. 
However, this also implies that the lower bound gets non-convex in general, which makes the fathoming tests significantly harder, and the lower bound improves only locally.

\subsection{Algorithmic Control of IP Scalarizations}\label{subsec:contr}
In the previous subsections we did not specify when to solve IP scalarizations, which implies a significant computational cost itself. However, this might be the most crucial part within the presented methods. Obviously, we aim at gaining as much information as possible by solving IP scalarizations. More objective space information will lead to tighter bounds that reduce the number of created nodes, due to a higher probability of fathoming by dominance and smaller search zones. Moreover, a reduced number of created nodes will reduce the total computation time. At the same time, solving overly many IP scalarizations will have a negative impact on the computation time. Furthermore, at a certain point the lower and upper bound will not improve anymore when solving additional IP scalarizations. 

So, there exists a trade-off between the reduction of the number of created subproblems and the decrease of the computation time. The difficulty is to find an appropriate condition to trigger an IP scalarization. Obviously, solving IP scalarizations more frequently in the beginning of the branch and bound algorithm is very promising. The earlier the lower and upper bounds are improved the more nodes might be fathomed. Moreover, solving the IP scalarization when the active node has weak bounds will lead to stronger improvements than in later stages of the algorithm. This is complemented by our adaptive branching strategy, which tends to select subproblems with weak lower bounds first.

The hybrid branch and bound algorithm using augmented weighted Tchebycheff scalarization entails also another problem. The augmented weighted Tchebycheff scalarization improves the lower bound just locally. If we use this scalarization at the beginning of the algorithm instead of the weighted sum scalarization, this could lead to an increase of created nodes. Once again, the intuitive idea is to start with the weighted sum IP scalarization more frequently in the beginning of the algorithm. This ensures that the lower bound improves globally at early stages of the branch and bound. The augmented weighted Tchebycheff scalarization should be used in later stages of the algorithm to find non-supported non-dominated points and to improve the lower bound locally. The efficiency of this idea and other approaches will be shown in the next section where we present numerical test results.

%% file: sec_res.tex
All algorithms were implemented in Julia 1.7.1 and the linear relaxations were solved with Bensolve 2.1 \citep{Loehne2017the}. The numerical tests were executed on a single core of a 3.20 GHz Intel\textsuperscript{\textregistered}  Core\texttrademark\ i7-8700 CPU processor in a computer with 32~GB RAM, running under openSUSE linux Leap~15.3.

We present numerical results of our new approaches and compare them to the general branch and bound framework presented in Section~\ref{sec:BB} which we use as baseline implementation. We consider three different types of problems: multidimensional knapsack problems, assignment problems and discrete facility location problems. The implementation of the proposed multiobjective branch and bound method and the considered benchmark instances are publicly available \citet{Bauss2023git}.
Multiple combinations of parameter settings are used to solve these test problems. Thereby, we compare the average number of explored nodes, the average number of solved IPs and the average computation time for 20 instances per problem size. The different evaluated approaches are 
\begin{itemize}
	\item the generic bi-objective Branch and Bouch (BB),
	\item bi-objective branch and bound using the local (BS1) respectively global (BS2) hypervolume gap as node selection criterion, 
	\item hybrid branch and bound including weighted sum IP scalarizations (WS), and 
	\item different combinations of the hybrid branch and bound algorithm using weighted sum IP scalarization (M1.$\alpha$.$\beta$) and hybrid branch and bound algorithm using weighted sum and augmented weighted Tchebycheff IP scalarization (M2.$\alpha$.$\beta$.$\gamma$).
\end{itemize}

The parameter $\alpha\in\{1,2,3\}$ controls how often IP scalarizations are applied. Since the number of IP scalarizations is chosen depending on the problem class, the meaning of the different values for \(\alpha\) is described in detail in the corresponding subsections. In general, however, the larger the parameter \(\alpha\) is chosen, the fewer IP scalarizations are solved.
With $\beta$ we distinguish between the local ($\beta =1$) and the global ($\beta=2$) hypervolume gap strategy. In the hybrid branch and bound algorithm using augmented weighted Tchebycheff scalarization we also distinguish between integrating the objective space information of the augmented weighted Tchebycheff into the lower bound ($\gamma=1$) or not ($\gamma=2$).

Note that the parameter values for each of the problem classes yield from preliminary test runs on a different sets of instances, where they shown to provide good results. Thus, the parameter values are chosen problem dependent but are not optimized for the specific test instances.

\subsection{Bi-objective Multidimensional Knapsack Problems}
We consider bi-objective, multidimensional knapsack problems with one, two and three linear restrictions (i.\,e.\ $m=1,2,3$). For every problem size we randomly generate 20 instances of the form
\[
\begin{array}{l r@{\extracolsep{0.75ex}}c@{\extracolsep{0.75ex}}l@{\extracolsep{2ex}}l}
	\max  & \displaystyle \sum_{i=1}^n c_i^k \, x_i   &&&\quad  k = 1,2\\
	\text{ s.t.}& \displaystyle \sum_{i=1}^n w_i\, x_i &\leq& b \\
	& \displaystyle \sum_{i=1}^n v_{ij} \,x_i &\leq& d_j &\quad j = 1,...,m-1\\
	& x &\in& \{0,1\}^n\\
\end{array}
\]
with $c^k_i \in[50,100]$, $w_i \in [ 5,15]$, $b = 5\, n, v_{ij}\in[ 5,15]$ and $d_j=\left\lfloor\frac{r\,n}{2}\right\rfloor$ with $r\in [5,15]$. 
Depending on the parameter \(\alpha\) we specify when and how often IP scalarizations are solved. In M1.1.$\beta$ and WS we apply the weighted sum scalarization every $10$-th iterations. In M1.2.$\beta$ we apply it every $10$-th iteration but only within the first $n^2$ iterations.  In M1.3.$\beta$ we apply the weighted sum scalarization every $10$-th iteration within the first $n^2/3$ iterations, every $n$-th iteration within the next $n^2/3$ iterations and every \mbox{$2n$-th} iteration within the third $n^2/3$ iterations. In M2.1.$\beta$.$\gamma$ we apply the weighted sum scalarization every $10$-th iteration and every $50$-th iteration the augmented weighted Tchebycheff scalarization is used instead. In M2.2.$\beta$.$\gamma$ we operate like in M1.2.$\beta$ but after the first $n^2$ iterations we apply the augmented weighted Tchebycheff scalarization every $50$-th iteration. In M2.3.$\beta$.$\gamma$ we operate like in M1.3.$\beta$ but after the first $n^2$ iterations we apply the augmented weighted Tchebycheff scalarization every $50$-th iteration. 
If a scalarization cannot be applied or the same IP scalarization has already been solved before, no IP scalarization is solved in that iteration.

\begin{table}
	\scriptsize
	\subfloat[Knapsack problem with $m=1$ constraint and $n=50$ variables	\label{tab1a}]{
		\begin{filecontents*}{Knapsack50.csv}
version,nodes,time (s),solved IPs
BB,27916.3,18.153,0.0
BS1,11788.1,8.339,0.0 
WS,14270.7,10.507,33.75
M1.1.1,10789.7,8.452,26.4 
M1.2.1,10793.5,8.188,21.2
M1.3.1,10795.7,8.116,17.95
M2.1.1.1,9888.5,10.873,48.7
M2.2.1.1,10140.3,9.437,32.65
M2.3.1.1,10521.0,8.774,25.65
M2.1.1.2,9840.1,8.396,45.55
M2.2.1.2,10130.6,8.422,32.35
M2.3.1.2,10401.8,8.288,26.25
BS2,16739.8,11.397,0.0 
M1.1.2,11026.3,8.861,26.1 
M1.2.2,11024.5,8.860,19.85
M1.3.2,11047.4,8.645,16.05
M2.1.2.1,10071.8,10.907,45.85
M2.2.2.1,10421.4,9.587,31.8
M2.3.2.1,10583.2,9.448,24.45
M2.1.2.2,9994.1,8.940,46.65
M2.2.2.2,10413.4,8.820,32.55
M2.3.2.2,10568.9,8.727,25.15
\end{filecontents*}
\pgfplotstabletypeset[
		col sep=comma,
		string type,
		every head row/.style={%
			before row={\hline
				\multicolumn{4}{|l|}{\rule{0pt}{1em}knapsack problem, $m=1,n=50$}  \\\hline
			},
			after row=\hline
		},
		every last row/.style={after row=\hline},
		columns/version/.style={column name=version, column type=|l},
		columns/nodes/.style={column name=nodes, column type=|C{1.2cm}},
		columns/time (s)/.style={column name=time (s), column type=|c},
		columns/solved IPs/.style={column name=\parbox{1cm}{IPs}, column type=|c|},
		]{Knapsack50.csv}
	}
	\hspace{0.25cm}
	\subfloat[Knapsack problem with $m=1$ constraint and $n=80$ variables\label{tab1b}]{
		\begin{filecontents*}{Knapsack80.csv}
version,nodes,time (s),solved IPs
BB,153938.9,186.330,0.0
BS1,36392.0,50.952,0.0 
WS,58825.7,79.545,54.0
M1.1.1,34337.7,50.431,41.65 
M1.2.1,34333.9,50.312,33.1
M1.3.1,34307.1,50.505,26.35
M2.1.1.1,31643.7,81.625,100.2
M2.2.1.1,32708.9,68.939,76.2
M2.3.1.1,32986.3,69.848,63.6
M2.1.1.2,31274.5,46.721,102.85
M2.2.1.2,32795.7,48.576,76.3
M2.3.1.2,33025.8,48.358,63.4
BS2,90976.0,116.847,0.0 
M1.1.2,39745.1,59.321,45.25 
M1.2.2,40083.2,59.350,31.2
M1.3.2,39918.1,58.999,24.5
M2.1.2.1,31905.8,80.505,99.7
M2.2.2.1,34496.9,79.444,84.0
M2.3.2.1,34571.7,72.955,65.15
M2.1.2.2,32074.9,48.510,104.85
M2.2.2.2,34169.8,51.464,87.15
M2.3.2.2,34943.3,51.887,63.1
\end{filecontents*}
\pgfplotstabletypeset[
		col sep=comma,
		string type,
		every head row/.style={%
			before row={\hline
				\multicolumn{4}{|l|}{\rule{0pt}{1em}knapsack problem, $m=1,n=80$}  \\\hline
			},
			after row=\hline
		},
		every last row/.style={after row=\hline},
		columns/version/.style={column name=version, column type=|l},
		columns/nodes/.style={column name=nodes, column type=|C{1.2cm}},
		columns/time (s)/.style={column name=time (s), column type=|c},
		columns/solved IPs/.style={column name=\parbox{1cm}{IPs}, column type=|c|},
		]{Knapsack80.csv}
	}

	\subfloat[Knapsack problem with $m=1$ constraint and $n=100$ variables \label{tab1c}]{
		\begin{filecontents*}{Knapsack100.csv}
version,nodes,time (s),solved IPs
BB,297345.3,484.676,0.0
BS1,68920.5,128.967,0.0 
WS,128080.8,224.587,66.95
M1.1.1,67369.1,128.665,54.1 
M1.2.1,67370.1,128.924,39.95
M1.3.1,67353.3,128.993,32.9
M2.1.1.1,58214.2,198.683,156.85
M2.2.1.1,61533.1,179.516,123.0
M2.3.1.1,62127.3,177.621,104.55
M2.1.1.2,58151.3,112.575,158.1
M2.2.1.2,61490.6,118.600,120.1
M2.3.1.2,61762.6,118.306,108.65
BS2,187306.9,318.524,0.0 
M1.1.2,73766.2,144.684,54.75 
M1.2.2,74065.4,144.677,37.9
M1.3.2,73865.0,144.306,31.0
M2.1.2.1,59512.7,200.754,158.05
M2.2.2.1,64489.2,192.211,127.0
M2.3.2.1,64330.5,187.803,114.65
M2.1.2.2,60470.8,118.479,157.75
M2.2.2.2,64943.8,127.428,123.5
M2.3.2.2,64711.1,126.525,113.5
\end{filecontents*}
\pgfplotstabletypeset[
		col sep=comma,
		string type,
		every head row/.style={%
			before row={\hline
				\multicolumn{4}{|l|}{\rule{0pt}{1em}knapsack problem, $m=1,n=100$}  \\\hline
			},
			after row=\hline
		},
		every last row/.style={after row=\hline},
		columns/version/.style={column name=version, column type=|l},
		columns/nodes/.style={column name= nodes, column type=|C{1.2cm}},
		columns/time (s)/.style={column name=time (s), column type=|c},
		columns/solved IPs/.style={column name=\parbox{1cm}{IPs}, column type=|c|},
		]{Knapsack100.csv}
	}
	\hspace{0.25cm}
	\subfloat[Knapsack problem with $m=2$ constraints and $n=50$ variables\label{tab1d}]{
		\begin{filecontents*}{BiKnapsack50.csv}
version,nodes,time (s),solved IPs
BB,32655.6,25.8684,0.0
BS1,10982.3,9.6578,0.0 
WS,14180.9,13.1749,33.25
M1.1.1,9784.7,9.5159,26.6 
M1.2.1,9782.5,9.2684,19.65
M1.3.1,9791.1,9.1580,14.85
M2.1.1.1,8900.7,12.6639,47.75
M2.2.1.1,9407.5,11.5112,34.5
M2.3.1.1,9507.5,11.0256,24.8
M2.1.1.2,8892.9,9.2702,47.0
M2.2.1.2,9370.2,9.2053,33.05
M2.3.1.2,9484.3,9.1161,24.75
BS2,15639.2,13.1246,0.0 
M1.1.2,10665.3,10.8423,28.5 
M1.2.2,10671.3,10.5141,19.4
M1.3.2,10854.2,10.5765,15.7
M2.1.2.1,9045.3,12.8916,48.9
M2.2.2.1,9629.7,12.1913,34.9
M2.3.2.1,9814.7,11.6068,28.1
M2.1.2.2,9030.0,9.5854,48.2
M2.2.2.2,9608.5,9.7621,36.05
M2.3.2.2,9787.6,9.6722,28.35
\end{filecontents*}
\pgfplotstabletypeset[
		col sep=comma,
		string type,
		every head row/.style={%
			before row={\hline
				\multicolumn{4}{|l|}{\rule{0pt}{1em}knapsack problem, $m=2,n=50$}  \\\hline
			},
			after row=\hline
		},
		every last row/.style={after row=\hline},
		columns/version/.style={column name=version, column type=|l},
		columns/nodes/.style={column name=nodes, column type=|C{1.2cm}},
		columns/time (s)/.style={column name=time (s), column type=|c},
		columns/solved IPs/.style={column name=\parbox{1cm}{IPs}, column type=|c|},
		]{BiKnapsack50.csv}
	}
	\caption{Numerical results of the bi-objective, multidimensional knapsack problems}
\end{table}		

\begin{table}
	\scriptsize
	\subfloat[Knapsack problem with $m=2$ constraints and $n=80$ variables\label{tab2a}]{
		\begin{filecontents*}{BiKnapsack80.csv}
version,nodes,time (s),solved IPs
BB,159911.4,287.925,0.0
BS1,41092.0,88.121,0.0 
WS,63215.0,130.338,55.0
M1.1.1,37799.1,82.654,43.9 
M1.2.1,37835.8,82.544,30.55
M1.3.1,37811.3,82.369,24.75
M2.1.1.1,31615.0,115.164,102.55
M2.2.1.1,34772.5,102.965,72.5
M2.3.1.1,35127.1,100.706,60.6
M2.1.1.2,31590.2,69.290,104.7
M2.2.1.2,34977.7,77.063,72.45
M2.3.1.2,35170.3,77.279,61.3
BS2,115223.3,224.926,0.0 
M1.1.2,43581.8,97.039,47.8 
M1.2.2,43744.8,97.874,29.1
M1.3.2,45481.1,102.689,23.45
M2.1.2.1,32388.8,116.173,106.25
M2.2.2.1,36453.2,120.161,78.3
M2.3.2.1,35942.7,116.972,69.05
M2.1.2.2,33264.8,74.207,104.75
M2.2.2.2,36578.1,81.915,77.2
M2.3.2.2,35505.1,77.971,69.55
\end{filecontents*}
\pgfplotstabletypeset[
		col sep=comma,
		string type,
		every head row/.style={%
			before row={\hline
				\multicolumn{4}{|l|}{\rule{0pt}{1em}knapsack problem, $m=2,n=80$}  \\\hline
			},
			after row=\hline
		},
		every last row/.style={after row=\hline},
		columns/version/.style={column name=version, column type=|l},
		columns/nodes/.style={column name= nodes, column type=|C{1.2cm}},
		columns/time (s)/.style={column name=time (s), column type=|c},
		columns/solved IPs/.style={column name=\parbox{1cm}{IPs}, column type=|c|},
		]{BiKnapsack80.csv}
	}
	\hspace{0.25cm}
	\subfloat[Knapsack problem with $m=2$ constraints and $n=100$ variables\label{tab2b}]{
		\begin{filecontents*}{BiKnapsack100.csv}
version,nodes,time (s),solved IPs
BB,428526.3,1074.21,0.0
BS1,100962.6,326.98,0.0 
WS,166108.1,464.71,67.25
M1.1.1,98831.5,323.54,54.8 
M1.2.1,99313.5,325.32,38.95
M1.3.1,98770.6,322.65,32.6
M2.1.1.1,69951.9,402.48,149.35
M2.2.1.1,73433.8,379.33,119.6
M2.3.1.1,73424.7,371.63,102.95
M2.1.1.2,70172.3,212.40,153.15
M2.2.1.2,72824.8,219.73,121.55
M2.3.1.2,73651.5,221.96,107.5
BS2,271110.8,720.46,0.0 
M1.1.2,113605.9,381.46,57.45 
M1.2.2,117188.3,394.57,36.45
M1.3.2,113665.3,378.63,29.85
M2.1.2.1,70603.0,404.28,150.8
M2.2.2.1,77836.0,400.18,121.4
M2.3.2.1,76818.2,399.89,110.8
M2.1.2.2,72316.9,219.91,148.4
M2.2.2.2,78135.2,240.57,121.3
M2.3.2.2,77073.0,235.49,112.45
\end{filecontents*}
\pgfplotstabletypeset[
		col sep=comma,
		string type,
		every head row/.style={%
			before row={\hline
				\multicolumn{4}{|l|}{\rule{0pt}{1em}knapsack problem, $m=2,n=100$}  \\\hline
			},
			after row=\hline
		},
		every last row/.style={after row=\hline},
		columns/version/.style={column name=version, column type=|l},
		columns/nodes/.style={column name= nodes, column type=|C{1.2cm}},
		columns/time (s)/.style={column name=time (s), column type=|c},
		columns/solved IPs/.style={column name=\parbox{1cm}{IPs}, column type=|c|},
		]{BiKnapsack100.csv}
	}

	\subfloat[Knapsack problem with $m=3$ constraints and $n=50$ variables\label{tab2c}]{
		\begin{filecontents*}{TriKnapsack50.csv}
version,nodes,time (s),solved IPs
BB,54430.4,51.5026,0.0
BS1,15260.1,17.9276,0.0 
WS,18112.9,20.5208,36.2
M1.1.1,13522.4,16.8431,28.8 
M1.2.1,13530.7,16.4735,19.0
M1.3.1,13576.5,16.4442,14.95
M2.1.1.1,12345.3,22.1289,51.15
M2.2.1.1,12973.5,19.7592,34.5
M2.3.1.1,13014.0,19.6348,28.8
M2.1.1.2,12241.1,16.1190,53.55
M2.2.1.2,12934.3,16.2933,35.15
M2.3.1.2,12908.3,16.1009,30.35
BS2,22597.3,24.5736,0.0 
M1.1.2,14645.7,18.6425,29.55 
M1.2.2,14573.2,18.0521,16.75
M1.3.2,14597.0,17.9793,14.5
M2.1.2.1,12617.9,22.3518,54.65
M2.2.2.1,13324.9,20.7655,33.4
M2.3.2.1,13252.3,20.4366,30.55
M2.1.2.2,12682.4,16.8670,56.65
M2.2.2.2,13180.2,16.7601,33.6
M2.3.2.2,13274.4,16.8497,32.15
\end{filecontents*}
\pgfplotstabletypeset[
		col sep=comma,
		string type,
		every head row/.style={%
			before row={\hline
				\multicolumn{4}{|l|}{\rule{0pt}{1em}knapsack problem, $m=3,n=50$}  \\\hline
			},
			after row=\hline
		},
		every last row/.style={after row=\hline},
		columns/version/.style={column name=version, column type=|l},
		columns/nodes/.style={column name= nodes, column type=|C{1.2cm}},
		columns/time (s)/.style={column name=time (s), column type=|c},
		columns/solved IPs/.style={column name=\parbox{1cm}{IPs}, column type=|c|},
		]{TriKnapsack50.csv}
	}
	\hspace{0.25cm}
	\subfloat[Knapsack problem with $m=3$ constraints and $n=80$ variables\label{tab2d}]{
		\begin{filecontents*}{TriKnapsack80.csv}
version,nodes,time (s),solved IPs
BB,263971.6,724.999,0.0
BS1,81609.9,287.899,0.0 
WS,121360.8,376.247,56.4
M1.1.1,80406.5,282.897,47.35 
M1.2.1,79971.5,279.885,32.2
M1.3.1,80089.6,279.686,26.55
M2.1.1.1,54187.1,340.389,115.7
M2.2.1.1,55915.9,328.411,85.45
M2.3.1.1,58486.8,330.347,66.9
M2.1.1.2,53396.0,164.755,116.0
M2.2.1.2,56572.6,174.131,86.25
M2.3.1.2,57452.9,176.657,67.85
BS2,140681.4,390.578,0.0 
M1.1.2,92175.8,334.414,50.45 
M1.2.2,96339.3,350.148,29.05
M1.3.2,96099.5,348.915,24.0
M2.1.2.1,54119.5,349.261,112.5
M2.2.2.1,59379.8,344.899,89.1
M2.3.2.1,60326.6,349.874,74.15
M2.1.2.2,54595.6,176.090,119.65
M2.2.2.2,62851.9,211.373,88.9
M2.3.2.2,61200.8,205.413,76.6
\end{filecontents*}
\pgfplotstabletypeset[
		col sep=comma,
		string type,
		every head row/.style={%
			before row={\hline
				\multicolumn{4}{|l|}{\rule{0pt}{1em}knapsack problem, $m=3,n=80$}  \\\hline
			},
			after row=\hline
		},
		every last row/.style={after row=\hline},
		columns/version/.style={column name=version, column type=|l},
		columns/nodes/.style={column name= nodes, column type=|C{1.2cm}},
		columns/time (s)/.style={column name=time (s), column type=|c},
		columns/solved IPs/.style={column name=\parbox{1cm}{IPs}, column type=|c|},
		]{TriKnapsack80.csv}
	}
	\caption{Numerical results of the bi-objective, multidimensional knapsack problems}
\end{table}

First of all, we notice that our branching strategies have a huge impact on the number of explored nodes and the computation time in knapsack problems. We observe that in general the local hypervolume gap strategy works better than the global hypervolume gap strategy. With the local strategy we can reduce the number of explored nodes by up to $76\%$ (Table~\ref{tab1c} and \ref{tab2b}) and the computation time by up to $73\%$ (Table~\ref{tab1c}). Although the local strategy works better the global hypervolume gap strategy has also a significant impact. The number of explored nodes can be reduced by up to $58\%$ (Table~\ref{tab2c}) and the computation time by up to $52\%$ (Table~\ref{tab2c}). The number of nodes and the computation time is reduced in all our approaches and we can notice that combinations with the local hypervolume strategy work better. 

By limiting the number of solved weighted sum IPs (i.\,e.\ in M1.2.$\beta$, M1.3.$\beta$, M2.2.$\beta$.$\gamma$ and M2.3.$\beta$.$\gamma$) we notice two consequences.  The number of nodes increases while the number of solved IPs decreases. Although  the number of nodes (and thus the number of considered subproblems) is increasing, the total computation time decreases. This implies that the reduced computation time to solve IP scalarizations compensates the increase of nodes, which results in a trade-off between the number of explored nodes and the computation time.  Another interesting aspect can be observed in M2.$\alpha$.$\beta$.1 and M2.$\alpha$.$\beta$.2. The computation time can be reduced if we do not integrate the augmented weighted Tchebycheff objective level set into the lower bound. This can be explained by the fact that the lower bound improvements of augmented weighted Tchebycheff are only local and do not compensate the computation time needed to integrate the information. The intuitive assumption that the number of explored nodes will then rise significantly is false. So, both our branching strategies work better, if we do not consider the local updates of the lower bound.

We can reach a reduction of the explored nodes by up to $83\%$ (Table~\ref{tab2b}) and a reduction of the computation time by up to $80\%$ (Table~\ref{tab2b}) in the best case. The strategies M2.1.1.1 and M2.1.1.2 seem to work best for knapsack problems. In most cases these two strategies have the largest impact on the number of explored nodes. Nevertheless, M2.1.1.2 achieves for all instance sizes the best computation times, since computation time is saved by not integrating the augmented weighted Tchebycheff objective space information into the lower bound.
Note that with  rising numbers of variables and constraints the hybridization techniques have larger impact on the performance of the branch and bound algorithm. 

\subsection{Bi-objective Assignment Problems}
We consider bi-objective assignment problems having $n=\ell^2$ variables,
\[
\begin{array}{l r@{\extracolsep{0.75ex}}c@{\extracolsep{0.75ex}}l@{\extracolsep{2ex}}l}
	\max  & \multicolumn{3}{l}{\displaystyle \sum_{i=1}^\ell \sum_{j=1}^\ell c_{ij}^k \, x_{ij}} & k =1,2 \\
	\text{ s.t.}& \displaystyle \sum_{i=1}^\ell x_{ij} &=& 1 & j = 1,...,\ell\\
	& \displaystyle  \sum_{j=1}^\ell x_{ij} &=& 1 & i = 1,...,\ell \\
	&x &\in& \{0,1\}^{\ell\times \ell} \\
\end{array}
\] 
where the cost coefficients $c_{ij}^k \in \left[50,100 \right]$.  The algorithmic strategy for the solution of IP scalarizations depending on the value of the parameter \(\alpha\) is chosen similarly to the previous case of knapsack problems. However, we adapt the boundaries due to the different number of nodes to explore in assignment problems.  In M1.1.$\beta$  weighted sum scalarizations are solved every $10$-th iteration to integer optimality. In M1.2.$\beta$ we apply the weighted sum every $10$-th iteration within the first $n\cdot\ell$ iterations.  In M1.3.$\beta$ we apply the weighted sum every $10$-th iteration within the first $n\cdot\ell/3$ iterations, every $\ell$-th iteration in the next  $n\cdot\ell/3$ iterations and every $n$-th iteration in the third $n\cdot\ell/3$ iterations.  For M.2.$\alpha$.$\beta$.$\gamma$ we use the same algorithmic strategy as in hybrid branch and bound for knapsack problems. 
If a scalarization cannot be applied or an IP with identical objective function has been solved prior, no IP is solved in that iteration.

\begin{table} 
	\scriptsize

	\subfloat[Assignment problem with $n=100$ variables\label{tab3a}]{
		\begin{filecontents*}{Assign10.csv}
version,nodes,time (s),solved IPs
BB,3117.0,5.1507,0.0
BS1,2422.2,4.1182,0.0 
WS,3117.0,5.2631,19.2
M1.1.1,2425.1,4.1937,13.95 
M1.2.1,2425.1,4.1564,11.95
M1.3.1,2425.1,4.1996,8.4
M2.1.1.1,2418.2,4.4063,19.8
M2.2.1.1,2420.5,4.2726,14.4
M2.3.1.1,2419.5,4.2242,9.65
M2.1.1.2,2420.5,4.3474,19.9
M2.2.1.2,2420.5,4.2014,14.55
M2.3.1.2,2423.0,4.2055,9.8
BS2,2948.9,4.9644,0.0 
M1.1.2,2519.8,4.4349,14.1 
M1.2.2,2518.7,4.4211,11.05
M1.3.2,2516.2,4.4203,7.65
M2.1.2.1,2499.7,4.6173,20.4
M2.2.2.1,2518.7,4.4661,12.2
M2.3.2.1,2516.2,4.4214,8.35
M2.1.2.2,2498.1,4.5400,19.55
M2.2.2.2,2518.7,4.4451,12.25
M2.3.2.2,2516.2,4.3972,8.4
\end{filecontents*}
\pgfplotstabletypeset[
		col sep=comma,
		string type,
		every head row/.style={%
			before row={\hline
				\multicolumn{4}{|l|}{\rule{0pt}{1em}assignment problem, $n=100$}  \\\hline
			},
			after row=\hline
		},
		every last row/.style={after row=\hline},
		columns/version/.style={column name=version, column type=|l},
		columns/nodes/.style={column name= nodes, column type=|C{1.2cm}},
		columns/time (s)/.style={column name=time (s), column type=|c},
		columns/solved IPs/.style={column name=\parbox{1cm}{IPs}, column type=|c|},,
		]{Assign10.csv}
	}
	\hspace{0.25cm}
	\subfloat[Assignment problem with $n=144$ variables\label{tab3b}]{
		\begin{filecontents*}{Assign12.csv}
version,nodes,time (s),solved IPs
BB,9661.9,28.2667,0.0
BS1,6255.4,18.9362,0.0 
WS,9610.8,28.1363,33.35
M1.1.1,6274.5,18.9323,22.3 
M1.2.1,6274.5,18.9869,14.3
M1.3.1,6274.5,18.9318,9.85
M2.1.1.1,6210.9,20.0216,46.0
M2.2.1.1,6279.9,19.4670,23.55
M2.3.1.1,6271.2,19.1542,12.65
M2.1.1.2,6171.0,19.4994,43.6
M2.2.1.2,6261.5,19.2851,24.65
M2.3.1.2,6270.5,18.8676,12.45
BS2,8448.0,24.7742,0.0 
M1.1.2,6781.8,20.6825,24.05 
M1.2.2,6779.4,20.5940,13.25
M1.3.2,6734.2,20.4270,9.3
M2.1.2.1,6472.8,20.8732,44.1
M2.2.2.1,6724.2,20.6838,16.8
M2.3.2.1,6682.3,20.4212,12.6
M2.1.2.2,6500.3,20.2536,42.55
M2.2.2.2,6737.7,20.5214,16.75
M2.3.2.2,6689.0,20.3987,12.95
\end{filecontents*}
\pgfplotstabletypeset[
		col sep=comma,
		string type,
		every head row/.style={%
			before row={\hline
				\multicolumn{4}{|l|}{\rule{0pt}{1em}assignment problem, $n=144$}  \\\hline
			},
			after row=\hline
		},
		every last row/.style={after row=\hline},
		columns/version/.style={column name=version, column type=|l},
		columns/nodes/.style={column name= nodes, column type=|C{1.2cm}},
		columns/time (s)/.style={column name=time (s), column type=|c},
		columns/solved IPs/.style={column name=\parbox{1cm}{IPs}, column type=|c|},
		]{Assign12.csv}
	}

	\subfloat[Assignment problem with $n=225$ variables\label{tab3c}]{
		\begin{filecontents*}{Assign15.csv}
version,nodes,time (s),solved IPs
BB,26810.5,160.712,0.0
BS1,16142.2,107.909,0.0 
WS,26810.5,163.666,46.85
M1.1.1,16150.0,107.842,32.8 
M1.2.1,16151.2,109.687,19.7
M1.3.1,16164.7,109.073,12.05
M2.1.1.1,15424.1,111.181,87.4
M2.2.1.1,15964.3,110.004,43.05
M2.3.1.1,16112.2,110.372,19.5
M2.1.1.2,15554.1,106.748,89.5
M2.2.1.2,16008.4,107.641,41.35
M2.3.1.2,16106.5,109.735,19.25
BS2,24566.2,152.341,0.0 
M1.1.2,17499.1,117.481,34.15 
M1.2.2,17591.4,118.804,16.75
M1.3.2,17287.8,116.681,11.3
M2.1.2.1,16095.7,115.012,85.8
M2.2.2.1,17048.8,116.275,31.55
M2.3.2.1,17093.3,117.976,17.75
M2.1.2.2,16183.0,110.920,79.55
M2.2.2.2,17104.7,117.885,32.2
M2.3.2.2,17070.0,117.557,16.65
\end{filecontents*}
\pgfplotstabletypeset[
		col sep=comma,
		string type,
		every head row/.style={%
			before row={\hline
				\multicolumn{4}{|l|}{\rule{0pt}{1em}assignment problem, $n=225$}  \\\hline
			},
			after row=\hline
		},
		every last row/.style={after row=\hline},
		columns/version/.style={column name=version, column type=|l},
		columns/nodes/.style={column name= nodes, column type=|C{1.2cm}},
		columns/time (s)/.style={column name=time (s), column type=|c},
		columns/solved IPs/.style={column name=\parbox{1cm}{IPs}, column type=|c|},
		]{Assign15.csv}
	}
	\hspace{0.25cm}
	\subfloat[Assignment problem with $n=324$ variables\label{tab3d}]{
		\begin{filecontents*}{Assign18.csv}
version,nodes,time (s),solved IPs
BB,76643.0,798.644,0.0
BS1,47311.6,527.471,0.0 
WS,75179.2,786.036,61.75
M1.1.1,47327.0,530.055,47.75 
M1.2.1,47332.8,523.562,21.9
M1.3.1,47375.6,524.610,15.25
M2.1.1.1,40978.0,477.927,157.5
M2.2.1.1,43760.5,497.812,82.0
M2.3.1.1,44343.9,499.678,52.35
M2.1.1.2,41294.4,457.120,159.05
M2.2.1.2,43864.6,487.051,81.25
M2.3.1.2,44499.7,489.744,52.05
BS2,69042.4,723.853,0.0 
M1.1.2,49367.5,565.719,48.45 
M1.2.2,49053.0,553.130,22.95
M1.3.2,49221.4,554.594,15.2
M2.1.2.1,41840.5,489.647,150.2
M2.2.2.1,45621.0,520.469,67.15
M2.3.2.1,46915.9,533.692,39.35
M2.1.2.2,42726.6,474.220,156.15
M2.2.2.2,45901.6,513.111,67.05
M2.3.2.2,46902.8,526.440,43.45
\end{filecontents*}
\pgfplotstabletypeset[
		col sep=comma,
		string type,
		every head row/.style={%
			before row={\hline
				\multicolumn{4}{|l|}{\rule{0pt}{1em}assignment problem, $n=324$}  \\\hline
			},
			after row=\hline
		},
		every last row/.style={after row=\hline},
		columns/version/.style={column name=version, column type=|l},
		columns/nodes/.style={column name= nodes, column type=|C{1.2cm}},
		columns/time (s)/.style={column name=time (s), column type=|c},
		columns/solved IPs/.style={column name=\parbox{1cm}{IPs}, column type=|c|},
		]{Assign18.csv}
	}
	\caption{Numerical results of the bi-objective assignment problems}
\end{table}

Due to the total unimodularity of the assignment problem, the weighted sum scalarizations do in general not improve the lower bound sets of subproblems. However, in situations where the weighted sum IP scalarization generates a supported efficient solution, whose corresponding non-dominated point is not an extreme point of the lower bound set, the local upper bounds move closer to the lower bound set. This reduces the gap between upper and lower bound and may lead to a decrease of the explored subproblems. Note that this update of the upper bound set may also have the contrary effect (the number of considered subproblems increases), since it can change the order in which the subproblems are considered. 
Although in general the weighted sum IP scalarizations are necessary to determine non-dominated points based on which the augmented weighted Tchebycheff scalarization can be applied, they are redundant for assignment problems due to the total unimodularity of their constraint matrix. 
Thus, all extreme supported non-dominated points (and their corresponding solutions) are obtained as extreme points of the lower bound set by the linear relaxation of the original problem which is solved in the root node of the branch and bound tree. However, we do not adjust our branch and bound algorithm for totally unimodular problem classes
to maintain comparable numerical results and general applicability.

Our branching strategies have a significant impact on the number of explored nodes and the computation time. Again, the local hypervolume gap performs better than the global hypervolume gap strategy. With the local strategy we can reduce the number of explored nodes by up to $39\%$ (Table~\ref{tab3c}) and the computation time by up to $33\%$ (Table~\ref{tab3c}). Using the global hypervolume gap strategy we can reduce the number of explored nodes by up to $12\%$ (Table~\ref{tab3b}) and the computation time by up to $12\%$ (Table~\ref{tab3b}). We reach a reduction of the explored nodes by up to $46\%$ (Table~\ref{tab3d}) and a reduction of the computation time by up to $42\%$ (Table~\ref{tab3d}), in the best case. Again, the strategies M2.1.1.1 and M2.1.1.2 seem to work the best for assignment problems in terms of explored nodes. Nevertheless, M2.1.1.2 leads to a better computation time which can be explained by the same argument as before. Furthermore, strategy BS1 works very well and is able to compete with the previously mentioned strategies with respect to number of nodes and computation time without solving a single IP scalarization.

\subsection{Bi-objective Discrete Facility Location Problems}
We consider discrete facility location problems of the form
\[
\begin{array}{l r@{\extracolsep{0.75ex}}c@{\extracolsep{0.75ex}}l@{\extracolsep{2ex}}l}
	\min   &\multicolumn{3}{l}{\displaystyle \sum_{i=1}^\ell \displaystyle\sum_{j=1}^q  c_{ij}^k\, x_{ij}  + \sum_{j=1}^q f_j^k \,y_j}  &\quad  k = 1,2\\
	\text{ s.t.}&  \displaystyle \sum_{j=1}^q x_{ij} &= & 1 & \quad i=1,...,\ell\\
	&x_{ij} &\leq& y_j &\quad i=1,...,\ell,\; j = 1,...,q\\[1.5ex]
	& x &\in &\{0,1\}^{\ell\times q}&\\[1.5ex]
	& y &\in & \{0,1\}^{q}&
\end{array}
\]
where $\ell$ is the number of customers and $q$ the number of facilities. We randomly generate coordinates of $\ell$ customers and $q$ facilities in a square with length $200$. The costs of the first objective function correspond to the $l_1$-distances between the customers and facilities, while the costs of the second objective function are randomly generated (i.\,e.\ $c_{ij}^2\in [1,200]$) and $f_j^k\in[200,400]$. The number of variables is $n=(\ell+1)\,q$. We restrict the numerical tests to problems where the number of facilities is $20\%$ of the number of customers. Again, we need to specify when and how often integer scalarizations are applied: We use the same methods as before but adapt the boundaries due to the different number of nodes to explore in discrete facility location problems. 
In M1.1.$\beta$ we apply the weighted sum IP scalarization every $10$-th iteration. In M1.2.$\beta$ we apply the weighted sum every $10$-th iteration within the first $n^2/4$ iterations.  In M1.3.$\beta$ we apply the weighted sum scalarization every $10$-th iteration in the first $n^2/4$ iterations, every $n/2$-th iteration in the next  $n^2/4$ iterations and every $n$-th iteration in the third $n^2/4$ iterations.  In M.2.$\alpha$.$\beta$.$\gamma$ we operate analogous to the methods used for knapsack and assignment problems. 
If a scalarization cannot be applied or an IP with identical objective function has been solved prior, no IP is solved in that iteration.

\begin{table}
	\scriptsize
	
	\subfloat[Facility location problem with $15$ customers and $3$ facilities\label{tab4a}]{
		\begin{filecontents*}{FacLoc15.3.csv}
version,nodes,time (s),solved IPs
BB,1431.1,1.0851,0.0
BS1,999.8,0.7640,0.0 
WS,1431.1,1.4550,15.35
M1.1.1,1000.2,0.8354,10.95 
M1.2.1,1000.2,0.8219,8.8
M1.3.1,1000.2,0.8243,6.95
M2.1.1.1,999.7,0.9361,13.2
M2.2.1.1,1000.2,0.8361,9.85
M2.3.1.1,1000.2,0.8079,7.85
M2.1.1.2,999.3,0.8536,12.8
M2.2.1.2,1000.2,0.8182,9.85
M2.3.1.2,1000.2,0.8037,7.85
BS2,1268.6,0.9714,0.0 
M1.1.2,1041.8,0.8654,12.3 
M1.2.2,1041.8,0.8527,9.35
M1.3.2,1041.8,0.8592,7.35
M2.1.2.1,1036.5,0.9722,15.65
M2.2.2.1,1041.8,0.8761,10.25
M2.3.2.1,1041.8,0.8735,7.6
M2.1.2.2,1034.6,0.9678,15.1
M2.2.2.2,1041.8,0.8843,10.4
M2.3.2.2,1041.8,0.8654,7.6
\end{filecontents*}
\pgfplotstabletypeset[
		col sep=comma,
		string type,
		every head row/.style={%
			before row={\hline
				\multicolumn{4}{|l|}{\rule{0pt}{1em}facility location problem, $n=48$}  \\\hline
			},
			after row=\hline
		},
		every last row/.style={after row=\hline},
		columns/version/.style={column name=version, column type=|l},
		columns/nodes/.style={column name= nodes, column type=|C{1.2cm}},
		columns/time (s)/.style={column name=time (s), column type=|c},
		columns/solved IPs/.style={column name=\parbox{1cm}{IPs}, column type=|c|},
		]{FacLoc15.3.csv}
	}
	\hspace{0.25cm}
	\subfloat[Facility location problem with $20$ customers and $4$ facilities\label{tab4b}]{
		\begin{filecontents*}{FacLoc20.4.csv}
version,nodes,time (s),solved IPs
BB,7949.1,12.6393,0.0
BS1,4626.7,7.9303,0.0 
WS,7490.2,12.2058,35.0
M1.1.1,4627.2,8.1249,26.45 
M1.2.1,4627.2,8.1149,18.2
M1.3.1,4626.4,8.0610,13.85
M2.1.1.1,4527.8,9.2394,49.0
M2.2.1.1,4601.6,8.4311,26.45
M2.3.1.1,4610.4,8.2569,18.9
M2.1.1.2,4526.1,8.4415,45.3
M2.2.1.2,4591.1,8.1289,24.25
M2.3.1.2,4605.2,8.1312,19.8
BS2,7084.4,11.5822,0.0 
M1.1.2,4873.8,8.7719,26.35 
M1.2.2,4874.8,8.7534,17.45
M1.3.2,4894.1,8.7031,12.6
M2.1.2.1,4673.6,9.5755,49.65
M2.2.2.1,4832.0,8.9275,23.6
M2.3.2.1,4812.8,8.8050,17.4
M2.1.2.2,4628.9,8.7019,46.15
M2.2.2.2,4825.9,8.7491,24.4
M2.3.2.2,4807.5,8.5984,17.5
\end{filecontents*}
\pgfplotstabletypeset[
		col sep=comma,
		string type,
		every head row/.style={%
			before row={\hline
				\multicolumn{4}{|l|}{\rule{0pt}{1em}facility location problem, $n=84$}  \\\hline
			},
			after row=\hline
		},
		every last row/.style={after row=\hline},
		columns/version/.style={column name=version, column type=|l},
		columns/nodes/.style={column name= nodes, column type=|C{1.2cm}},
		columns/time (s)/.style={column name=time (s), column type=|c},
		columns/solved IPs/.style={column name=\parbox{1cm}{IPs}, column type=|c|},
		]{FacLoc20.4.csv}
	}

	\subfloat[Facility location problem with $25$ customers and $5$ facilites\label{tab4c}]{
		\begin{filecontents*}{FacLoc25.5.csv}
version,nodes,time (s),solved IPs
BB,17461.9,51.6307,0.0
BS1,10684.0,34.5157,0.0 
WS,16795.9,50.9317,51.35
M1.1.1,10753.9,35.4356,37.1 
M1.2.1,10753.9,35.2764,25.5
M1.3.1,10753.9,35.1007,19.15
M2.1.1.1,10104.4,38.3909,89.65
M2.2.1.1,10678.1,35.7434,39.35
M2.3.1.1,10722.3,35.6262,27.7
M2.1.1.2,10103.0,34.5003,85.95
M2.2.1.2,10691.2,35.3730,39.05
M2.3.1.2,10718.6,35.1690,27.45
BS2,15474.7,46.5130,0.0 
M1.1.2,11548.8,38.4891,39.3 
M1.2.2,11601.4,38.5848,24.6
M1.3.2,11535.1,38.2917,18.0
M2.1.2.1,10684.3,39.8639,81.5
M2.2.2.1,11381.8,39.1988,37.15
M2.3.2.1,11336.0,38.7754,28.45
M2.1.2.2,10695.5,36.9098,78.25
M2.2.2.2,11341.6,38.2646,39.55
M2.3.2.2,11352.1,38.0063,25.9
\end{filecontents*}
\pgfplotstabletypeset[
		col sep=comma,
		string type,
		every head row/.style={%
			before row={\hline
				\multicolumn{4}{|l|}{\rule{0pt}{1em}facility location problem, $n=130$}  \\\hline
			},
			after row=\hline
		},
		every last row/.style={after row=\hline},
		columns/version/.style={column name=version, column type=|l},
		columns/nodes/.style={column name= nodes, column type=|C{1.2cm}},
		columns/time (s)/.style={column name=time (s), column type=|c},
		columns/solved IPs/.style={column name=\parbox{1cm}{IPs}, column type=|c|},
		]{FacLoc25.5.csv}
	}
	\hspace{0.25cm}
	\subfloat[Facility location problem with $30$ customers and $6$ facilities\label{tab4d}]{
		\begin{filecontents*}{FacLoc30.6.csv}
version,nodes,time (s),solved IPs
BB,67369.3,373.238,0.0
BS1,31844.1,203.145,0.0 
WS,62192.8,349.741,69.95
M1.1.1,32106.6,206.244,53.05 
M1.2.1,32097.9,206.851,33.4
M1.3.1,32148.5,207.199,23.75
M2.1.1.1,28384.2,224.486,186.8
M2.2.1.1,31074.5,218.314,101.15
M2.3.1.1,32004.8,209.608,50.1
M2.1.1.2,28558.8,186.010,172.7
M2.2.1.2,30946.4,202.329,97.75
M2.3.1.2,32011.8,207.715,47.5
BS2,50759.0,292.344,0.0 
M1.1.2,35150.1,230.687,55.5 
M1.2.2,35789.8,233.967,30.8
M1.3.2,35255.0,233.163,21.95
M2.1.2.1,29704.3,230.016,172.65
M2.2.2.1,33959.7,236.351,81.75
M2.3.2.1,34228.0,233.553,50.2
M2.1.2.2,30412.6,202.719,167.4
M2.2.2.2,34511.8,229.970,69.95
M2.3.2.2,34405.0,229.021,47.1
\end{filecontents*}
\pgfplotstabletypeset[
		col sep=comma,
		string type,
		every head row/.style={%
			before row={\hline
				\multicolumn{4}{|l|}{\rule{0pt}{1em}facility location problem, $n=186$}  \\\hline
			},
			after row=\hline
		},
		every last row/.style={after row=\hline},
		columns/version/.style={column name=version, column type=|l},
		columns/nodes/.style={column name= nodes, column type=|C{1.2cm}},
		columns/time (s)/.style={column name=time (s), column type=|c},
		columns/solved IPs/.style={column name=\parbox{1cm}{IPs}, column type=|c|},
		]{FacLoc30.6.csv}
	}
	\caption{Numerical results of the bi-objective integer facility location problem}
\end{table}

Again, both new branching strategies have an impact on the number of explored nodes and the computation time. The local strategy, once more, leads to better results, namely reduction of the explored nodes by up to $52\%$ (Table~\ref{tab4d}) and reduction of the computation time by up to $45\%$ (Table~\ref{tab4d}). With the global hypervolume gap strategy we can reach a reduction of the explored nodes by up to $24\%$ (Table~\ref{tab4d}) and a reduction of the computation time by up to $21\%$ (Table~\ref{tab4d}). In the best case we can reach a reduction of the explored nodes by up to $57\%$ (Table~\ref{tab4d}) and of the computation time by up to $50\%$ (Table~\ref{tab4d}). 
Once again, M2.1.1.2 seems to be the best choice with respect to the number of explored nodes and with a rising number of variables it is also the best choice regarding the computation time. With a smaller number of variables, BS1 leads to good results with respect to both aspects without solving a single IP.

\subsection{Summary}
In all of the three tested problem classes (knapsack, assignment, discrete facility location) a significant reduction of the number of explored nodes and the computation time can be realized with all presented combinations of the hybrid branch and bound approach. 
With increasing problem size (number of variables) the impact of the presented augmentations increases. Furthermore, the approaches perform better on problems where the gap between $Y_N$ and the solution of the linear relaxation is larger compared to \emph{totally unimodular} problems. The reduction in terms of the number of branch and bound nodes and runtime we achieve with the proposed methods as compared to plain branch and bound is visualized in Figure \ref{fig:plots} for varying instance sizes.

The local hypervolume gap strategy for the node selection outperforms the global hypervolume gap strategy in our numerical tests. A reason for this is that in the global hypervolume gap strategy many small search zones can add up to a large gap although the lower bound might be quite close to the non-dominated points. The local hypervolume gap strategy chooses the node with the largest search zone, which has the biggest potential to reduce this gap. Moreover, the local hypervolume gap strategy aims at an uniform distribution of points in the incumbent list.
In our numerical test, M2.1.1.2 turn out to be the best choice in most cases with respect to the number of explored nodes and computation time. In this version, we use the local hypervolume gap strategy for the choice of the active node, every $10$-th iteration the weighted sum IP scalarization is applied and every $50$-th iteration we apply the augmented weighted Tchebycheff scalarization instead. Futhermore, the objective space information gained by the augmented weighted Tchebycheff scalarization is not used to update the lower bound set, since its local improvements do not compensate the increased computation time. Although we need to solve more IPs than in most other approaches, the computation time is the lowest compared to the others. So, using the augmented weighted Tchebycheff scalarization in the beginning of the branch and bound works best. Due to the likelihood of finding non-supported non-dominated points in the early stages of the algorithm, the upper bound can be further improved. This results to a higher probability of fathoming a node by dominance.  Nevertheless, with version BS1 we also achive a remarkable reduction in terms of the number of explored nodes and computation time by using the local hypervolume gap strategy for node selection.

	\begin{center}
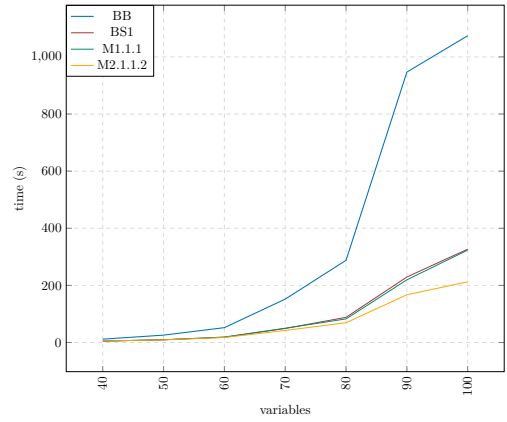
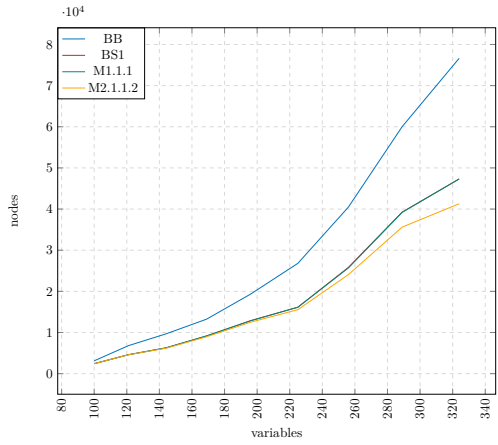
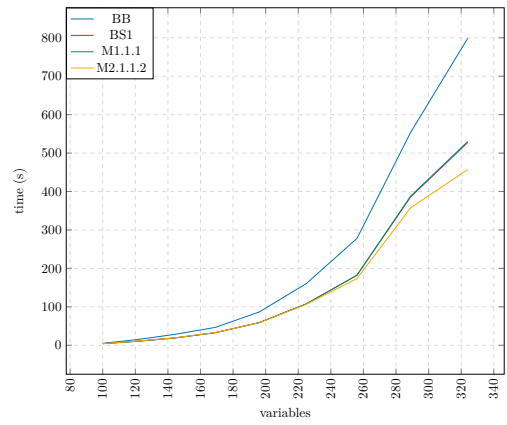
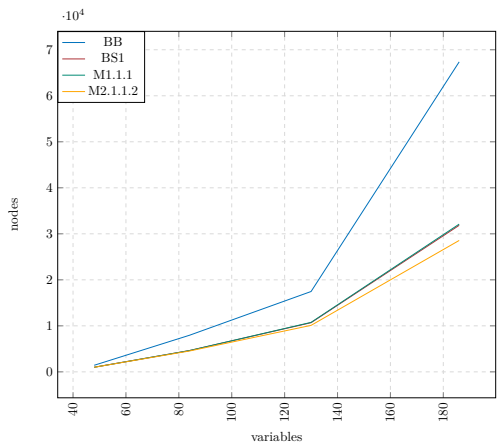
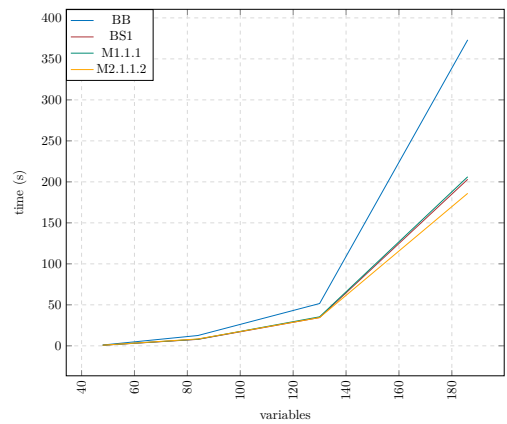
\begin{figure}[htbp!]
	\subfloat[Number of nodes for 2-D knapsack problems\label{FigBiKnaNodes}]{
		\begin{tikzpicture}[scale =0.499]
			\begin{axis}[
				width=\linewidth, 
				grid=major, 
				grid style={dashed,gray!30},
				y label style={at={(axis description cs:-0.08,.5)},anchor=south},
				xlabel=variables, 
				ylabel=nodes,
				legend style={at={(0,1)},anchor=north west},
				x tick label style={rotate=90,anchor=east}
				]
				\addplot[color = RoyalBlue]
				table[x=variables,y=BB,col sep=comma] {BiKnapPlotNodes.csv}; 
				\addplot[color = Maroon]
				table[x=variables,y=BS1,col sep=comma] {BiKnapPlotNodes.csv}; 
					\addplot[color = PineGreen]
				table[x=variables,y=M111,col sep=comma] {BiKnapPlotNodes.csv}; 
				\addplot[color = ChromeYellow]
				table[x=variables,y=M2112,col sep=comma] {BiKnapPlotNodes.csv}; 
				egend style={anchor=north east}
				\legend{BB,BS1,M1.1.1,M2.1.1.2}
				
			\end{axis}
		\end{tikzpicture}
	}
\hspace{0.25cm}
	\subfloat[Runtime for 2-D knapsack problems\label{FigBiKnaTime}]{
	\begin{tikzpicture}[scale =0.499]
		\begin{axis}[
			width=\linewidth, 
			grid=major, 
			grid style={dashed,gray!30},
			y label style={at={(axis description cs:-0.08,.5)},anchor=south},
			xlabel=variables, 
			ylabel= time (s),
			legend style={at={(0,1)},anchor=north west},
			x tick label style={rotate=90,anchor=east}
			]
			\addplot[color = RoyalBlue,thick]
			table[x=variables,y=BB,col sep=comma] {BiKnapPlotTime.csv}; 
			\addplot[color = Maroon]
			table[x=variables,y=BS1,col sep=comma] {BiKnapPlotTime.csv}; 
			\addplot[color = PineGreen]
			table[x=variables,y=M111,col sep=comma] {BiKnapPlotTime.csv}; 
			\addplot[color = ChromeYellow]
			table[x=variables,y=M2112,col sep=comma] {BiKnapPlotTime.csv}; 
			
			\legend{BB,BS1,M1.1.1,M2.1.1.2}
			
		\end{axis}
	\end{tikzpicture}
}

	\subfloat[Number of nodes for assignment problems\label{FigAssignNodes}]{
	\begin{tikzpicture}[scale =0.499]
		\begin{axis}[
			width=\linewidth, 
			grid=major, 
			grid style={dashed,gray!30},
			y label style={at={(axis description cs:-0.08,.5)},anchor=south},
			xlabel=variables, 
			ylabel= nodes,
			legend style={at={(0,1)},anchor=north west},
			x tick label style={rotate=90,anchor=east}
			]
			\addplot[color = RoyalBlue]
			table[x=variables,y=BB,col sep=comma] {AssignPlotNodes.csv}; 
			\addplot[color = Maroon]
			table[x=variables,y=BS1,col sep=comma] {AssignPlotNodes.csv}; 
			\addplot[color = PineGreen]
			table[x=variables,y=M111,col sep=comma] {AssignPlotNodes.csv}; 
			\addplot[color = ChromeYellow]
			table[x=variables,y=M2112,col sep=comma] {AssignPlotNodes.csv}; 
			
			\legend{BB,BS1,M1.1.1,M2.1.1.2}
			
		\end{axis}
	\end{tikzpicture}
}
	\hspace{0.25cm}
\subfloat[Runtime for assignment problems\label{FigAssignTime}]{
	\begin{tikzpicture}[scale =0.499]
		\begin{axis}[
			width=\linewidth, 
			grid=major, 
			grid style={dashed,gray!30},
			y label style={at={(axis description cs:-0.08,.5)},anchor=south},
			xlabel=variables, 
			ylabel= time (s),
			legend style={at={(0,1)},anchor=north west},
			x tick label style={rotate=90,anchor=east}
			]
			\addplot[color = RoyalBlue]
			table[x=variables,y=BB,col sep=comma] {AssignPlotTime.csv}; 
			\addplot[color = Maroon]
			table[x=variables,y=BS1,col sep=comma] {AssignPlotTime.csv}; 
			\addplot[color = PineGreen]
			table[x=variables,y=M111,col sep=comma] {AssignPlotTime.csv}; 
			\addplot[color = ChromeYellow]
			table[x=variables,y=M2112,col sep=comma] {AssignPlotTime.csv}; 
			
			\legend{BB,BS1,M1.1.1,M2.1.1.2}
			
		\end{axis}
	\end{tikzpicture}
}

\subfloat[Number of nodes for facility location problems\label{FigFLNodes}]{
	\begin{tikzpicture}[scale =0.499]
		\begin{axis}[
			width=\linewidth, 
			grid=major, 
			grid style={dashed,gray!30},
			y label style={at={(axis description cs:-0.08,.5)},anchor=south},
			xlabel=variables, 
			ylabel= nodes,
			legend style={at={(0,1)},anchor=north west},
			x tick label style={rotate=90,anchor=east}
			]
			\addplot[color = RoyalBlue]
			table[x=variables,y=BB,col sep=comma] {FLPlotNodes.csv}; 
			\addplot[color = Maroon]
			table[x=variables,y=BS1,col sep=comma] {FLPlotNodes.csv}; 
			\addplot[color = PineGreen]
			table[x=variables,y=M111,col sep=comma] {FLPlotNodes.csv}; 
			\addplot[color = ChromeYellow]
			table[x=variables,y=M2112,col sep=comma] {FLPlotNodes.csv}; 
			
			\legend{BB,BS1,M1.1.1,M2.1.1.2}
			
		\end{axis}
	\end{tikzpicture}
}
	\hspace{0.25cm}
\subfloat[Runtime for facility location problems\label{FigFLTime}]{
	\begin{tikzpicture}[scale =0.499]
		\begin{axis}[
			width=\linewidth, 
			grid=major, 
			grid style={dashed,gray!30},
			y label style={at={(axis description cs:-0.08,.5)},anchor=south},
			xlabel=variables, 
			ylabel= time (s),
			legend style={at={(0,1)},anchor=north west},
			x tick label style={rotate=90,anchor=east}
			]
			\addplot[color = RoyalBlue]
			table[x=variables,y=BB,col sep=comma] {FLPlotTime.csv}; 
			\addplot[color = Maroon]
			table[x=variables,y=BS1,col sep=comma] {FLPlotTime.csv}; 
			\addplot[color = PineGreen]
			table[x=variables,y=M111,col sep=comma] {FLPlotTime.csv}; 
			\addplot[color = ChromeYellow]
			table[x=variables,y=M2112,col sep=comma] {FLPlotTime.csv}; 
			
			\legend{BB,BS1,M1.1.1,M2.1.1.2}
			
		\end{axis}
	\end{tikzpicture}
}
\caption{Visualization of branch and bound node reduction and runtime reduction for varying test instance sizes on a selection of approaches.}
\label{fig:plots}
\end{figure}

	\end{center}

%% file: sec_concl.tex
In this paper, we propose two approaches to incorporate objective space information in bi-objective branch and bound. By using the local or global (approximated) hypervolume gap as a node selection criterion, we adapt the run of the branch and bound algorithm to the problem instance. Additionally, we adaptively solve scalarizations to integer optimality to improve the lower and the upper bound set by the obtained objective space information. Our numerical results show the effectiveness of both approaches and in particular of their combination. The dynamic branching rule based on the local (approximated) hypervolume gap has large impact on the number of explored subproblems, is compuationally efficient and can be easily integrated in other multi-objective branch and bound algorithms.

While we tested in this paper the individual contributions of our augmentations on a generic bi-objective branch and bound, we will continue to extend our ideas to multiple dimensions and integrate them into a competetive multi-objective branch and bound framework. Particularly in higher dimensions, it may be promising to combine our approaches with objective space branching.